\documentclass[a4paper,10pt]{article}
\usepackage[utf8x]{inputenc}
\usepackage{latexsym}
\usepackage{epsfig}
\usepackage{amsmath}
\usepackage{amssymb}
\usepackage{amsfonts}
\usepackage{amsbsy}
\usepackage{graphicx}
\usepackage{graphics}
\usepackage{euscript}
\usepackage{epsfig}
\usepackage{varioref}
\usepackage{enumerate}
\usepackage{verbatim}
\usepackage{makeidx}
\usepackage{color}

\makeindex

\def\R{\mathbb{R}}
\def\T{\mathbb{T}}
\def\C{\mathbb{C}}

\def\P{\mathbb{P}}

\def\N{\mathbb{N}}
\def\u{{\bf u}}
\def\v{{\bf v}}
\def\w{{\bf w}}

\def\n{{\bf n}}
\def\x{{\bf x}}
\def\f{{\bf f}}

\newtheorem{theorem}{Theorem}[section]
\newtheorem{proposition}{Proposition}[section]

\newtheorem{remark}{Remark}[section]
\newtheorem{corollary}{Corollary}[section]
\newtheorem{lemma}{Lemma}[section]
\newtheorem{assumption}{Assumption}[section]
\numberwithin{equation}{section}

\makeatletter
\renewcommand\section{\@startsection {section}{1}{\z@}%
 				   {-2.5ex \@plus -1ex \@minus -.2ex}%
                                   {1.3ex \@plus.2ex}%
                                   {\normalfont\large\bfseries}}
\renewcommand\subsection{\@startsection {subsection}{1}{\z@}%
 				   {-2.ex \@plus -1ex \@minus -.2ex}%
                                   {0.5ex \@plus.2ex}%
                                   {\normalfont\normalsize\bfseries}}
\makeatother

\begin{document}

\title{Stability conditions for the numerical solution of convection-dominated problems with skew-symmetric discretizations}

\author{
Erwan Deriaz\thanks{Laboratoire de M\'ecanique, Mod\'elisation et Proc\'ed\'es Propres, 38 rue Fr\'ed\'eric Joliot-Curie
13451 MARSEILLE Cedex 20 (France), erwan.deriaz@l3m.univ-mrs.fr}
}


\bigskip
\bigskip

\date{\today}

\maketitle

\begin{abstract}
This paper presents original and close to optimal stability conditions linking the time step and the space step,
stronger than the CFL criterion: $\delta t\leq C\delta x^\alpha$ with $\alpha=\frac{2r}{2r-1}$, $r$ an integer,
for some numerical schemes we produce, when solving convection-dominated problems. We test this condition
numerically and prove that it applies to nonlinear equations under smoothness assumptions.

\end{abstract}

{\bf keywords:} CFL condition, von Neumann stability, transport equation, Euler equation, Runge-Kutta schemes, Adams-Bashforth schemes.

\section{Introduction}
\label{intro}

In numerical fluid mechanics, many simulations for transport-dominated problems employ explicit second order time discretization schemes,
either of Runge-Kutta type \cite{KT97,DP07} or Adams-Bashforth \cite{Pey00,Sch05}.
Although widely in use and proved efficient, the stability domains of these order two numerical schemes (see Fig. \ref{RK&AB})
exclude the $(Oy)$ axis corresponding to transport problems.
Nonetheless, actual experiments \cite{ZS04,DP07} show that even in this case, a convergent solution can be obtained.
If the problem admits a sufficiently smooth, classical solution, the second order time-stepping is stable
at worst under a condition
of type $\delta t\leq C(\delta x/u_{\rm max})^{4/3}$,
where $\delta t$ is the time step, $\delta x$ the space step, and $u_{\rm max}$ the maximum velocity of the transport problem.

A close look at the stability condition provided by an analysis of von Neumann type applied to transport equation provides an explanation.
To the best of our knowledge, this result is new --for instance, it is not presented in \cite{Wes01} which has collected
the state of the art in numerical stability-- despite
the fact that it applies to a wide variety of numerical problems.
Under some smoothness conditions, it readily extends to Burgers equation, incompressible Euler equations, Navier-Stokes
equations with a high Reynolds number on domains possibly bounded by walls, and
to conservation laws.

For the single step numerical method (i.e. the explicit Euler scheme), a stability result relying on a similar approach
and providing a stability constraint of the type $\delta t\leq C(\delta x/u_{\rm max})^{2}$ has been presented by
several authors \cite{GT91,MT98,Tre94}.
The square originates from a completely different kind of numerical instability than the usual
stability condition for the heat equation with explicit schemes.
As we will see in this article, it comes from the order of tangency of the stability domain to the $(Oy)$
axis and applies only to some first order schemes while for the heat equation it comes from the second derivative
notwithstanding the order of the scheme.
We present the generalization of this stability constraint to other schemes.
Incidentally, we show that for transport dominated problems there exists a direct connection between the order
and the stability of numerical schemes.

As the numerical viscosity may stabilize the time scheme, this $\frac{4}{3}$-CFL\footnote{CFL stands for the names of the three authors of the
founding paper \cite{CFL67}: R.~Courant, K.~Friedrichs and H.~Lewy} criterion applies essentially
to pseudo-spectral methods and conservative numerical methods \cite{VV03,ZS04}.
A basic numerical experiment allows us to validate our approach.

The paper is organized as follows: first we recall the definition and the computation of the
von Neumann stability; then we focus on the linear transport problem, predicting a stability
condition of the type $\delta t \leq C(\delta x/u)^{2r/(2r-1)}$ with $r$ an integer, for several schemes; then we construct
numerical schemes for which such a stability condition appears for $r=1,2,3,4$, and corresponds to exponents equal to $2$, $\frac{4}{3}$,
$\frac{6}{5}$ and $\frac{8}{7}$; finally we show how this stability criterion extends to nonlinear equations,
and to multicomponent transport equations (including wave equations).

\section{The von Neumann stability condition}
\label{VNstab}
Let us consider the equation
\begin{equation}
\label{diffequ}
\partial_t u = F\,u,\quad u(0,\cdot)=u_0
\end{equation}
where $u:\R_+\times\R\to\R,\,(t,x)\mapsto u(t,x)$ and $F$ is a linear operator. We denote by $\sigma(\xi)$ the symbol
associated to $F$, i.e. $\widehat{F\,u}(t,\xi)=\sigma(\xi)\widehat{u}(t,\xi)$ where $u\mapsto \widehat{u}$
stands for the Fourier transform\footnote{The Fourier transform of a
function $f\in L^1(\R)$ is noted $\hat f (\xi)=\int_{-\infty}^{+\infty} f(x) ~e^{-ix\xi}
dx$, we recall that $f\mapsto \frac{1}{\sqrt{2\pi}}\, \hat f$
defines an isometry on $L^2(\R)$.}.

\ 

In the following, we explain how to apply the von Neumann stability analysis as presented in \cite{Pey00,Tre94,Wes01}.
We note $u_k\sim u(k\delta t,\cdot)$ the approximation at time $k\delta t$
for $k$ an integer, $\delta t$ denoting the time step.
We consider we have a spectral discretization,
or that all the terms are orthogonally reprojected in our discretization space.
The scheme can be of Runge-Kutta type, relying on the computation of intermediate time steps $u_{(\ell)}$:
\begin{equation}
\label{RKtype}
u_{(0)}=u_n,\quad
u_{(\ell)}=\sum_{i=0}^{\ell-1}a_{\ell i}\,u_{(i)}+\delta t \,\sum_{i=0}^{\ell-1}b_{\ell i}\,F\,u_{(i)}
~~{\rm for}~1\leq \ell \leq s',\quad
u_{n+1}=u_{(s')}
\end{equation}
with $(a_{\ell i})_{\ell,i}$ and $(b_{\ell i})_{\ell,i}$ well chosen to ensure
the accuracy of the integration.

Or it can be an explicit multi-step (Adams-Bashforth) scheme involving the previous time steps:
\begin{equation}
\label{ABtype}
u_{n+1}=\sum_{i=0}^{s}c_{i}\,u_{n-i}+\delta t \,\sum_{i=0}^{s}d_{i}\,F\,u_{n-i}.
\end{equation}
We can also mix these two types of integration schemes:
\begin{eqnarray}
\nonumber
u_{(\ell)}=\sum_{i=1}^{\ell-1}a_{\ell i}\,u_{(i)}+\sum_{i=0}^{s}c_{\ell i}\,u_{n-i}+\delta t \,\sum_{i=0}^{\ell-1}b_{\ell i}\,F\,u_{(i)}
+\delta t \,\sum_{i=0}^{s}d_{\ell i}\,F\,u_{n-i}
~{\rm for}~1\leq \ell \leq s',\\
\label{mixtype} u_{n+1}=u_{(s')}.~~~~~~~~~~~
\end{eqnarray}
The von Neumann stability analysis consists in isolating a Fourier mode $\xi$ by taking $u_n(x)=\phi_n\,e^{i\xi x}$.
Actually, if $\delta x$ is the space step, then $-\frac{\pi}{\delta x}\leq \xi \leq\frac{\pi}{\delta x}$.\\
In the case when several previous time samples are necessary, like in the case of an Adams-Bashforth scheme, we set
\begin{equation}
\label{X_n}
X_n=\left( \begin{array}{c}
u_n \\ u_{n-1} \\ \vdots \\ u_{n-s}
\end{array} \right) .
\end{equation}
Remarking that each time we apply $F$ to a term in (\ref{mixtype}), we also multiply this term by $\delta t$,
it turns out that
\begin{equation}
\label{defMsigma}
X_{n+1}=M(\sigma(\xi)\delta t)X_n
\end{equation}
where, setting $\zeta=\sigma(\xi)\delta t$, $M(\zeta)$ is a $(s+1)\times(s+1)$ square matrix
whose elements are polynomials in $\zeta$.
Note that if $F$ is a differential operator with derivatives of maximal order $\gamma$, then $|\zeta|\leq K \frac{\delta t}{\delta x^\gamma}$.
In the case of hyperbolic equations, $\gamma$ is equal to one.

Let $\lambda_0(\zeta),\dots,\lambda_{s}(\zeta)$ denote the eigenvalues of $M(\zeta)$.
The spectral radius is defined by \begin{equation}\rho(M(\zeta))={\rm max}_{0\leq i \leq s}|\lambda_i(\zeta)|.\end{equation}
Then \begin{equation}\rho(M(\zeta))^n\leq \|M(\zeta)^n \|\leq \|M(\zeta) \|^n.\end{equation}
For almost every $\zeta$, $\exists K_\zeta>0$ such that
$\forall n\geq 0$, $\|M(\zeta)^n \|\leq K_\zeta\,\rho(M(\zeta))^n$ where
the constant $K_\zeta$ becomes large near the singularities of $M(\zeta)$.
Actually, in numerical experiments, this does not play any crucial role (see \cite{Tre94} for a complete discussion on this topic).
Hence, overlooking this latest point, the von Neumann stability of the scheme (\ref{mixtype}) is assured by:
\begin{equation}
\label{condMxi}
\forall i,\zeta, \quad |\lambda_i(\zeta)|\leq 1 +C\delta t
\end{equation}
with $C$ a positive constant independent of $\delta x$ and $\delta t$. Sometimes, $C$ is taken equal to zero to enforce
an absolute stability.
The assumption (\ref{condMxi}) allows any error $\varepsilon_0$ to stay bounded after an elapsed time $T$, since:
\begin{equation}
\|\varepsilon_T\|=\| M(\zeta)^{T/\delta t} \varepsilon_0\| \leq K_\zeta\,(1 +C\delta t)^{T/\delta t} \|\varepsilon_0\|
\leq K_\zeta\,e^{CT} \|\varepsilon_0\|.
\end{equation}

\ 

\begin{figure}
\begin{center}
\begin{tabular}{cc}
\includegraphics[scale=0.42]{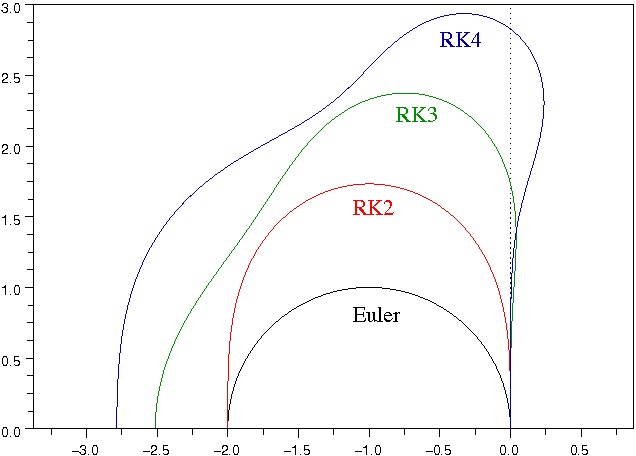} &
\includegraphics[scale=0.42]{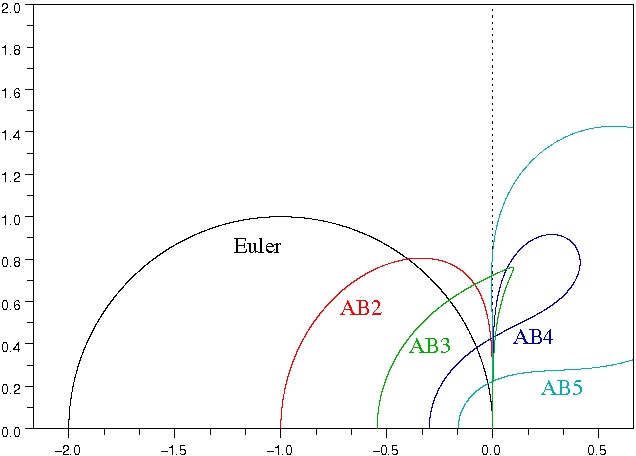} \\
Runge-Kutta & Adams-Bashforth
\end{tabular}
\caption{\label{RK&AB} Von Neumann stability domains for the first four Runge-Kutta and five Adams-Bashforth schemes.}
\end{center}
\end{figure}

\begin{figure}
\centerline{\includegraphics[width=10cm,clip]{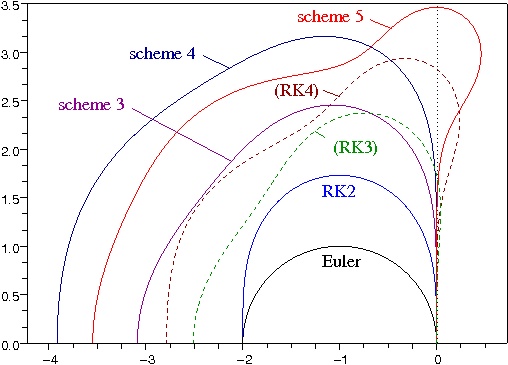}}
\caption{\label{neu_sch2N} Von Neumann stability domains for schemes of Runge-Kutta type.}
\end{figure}

The von Neumann {\it stability domain} of the scheme (\ref{mixtype}) is given by $\mathcal{S}=\{\zeta\in\C,\rho(M(\zeta))\leq 1\}$.
In Fig. \ref{RK&AB}, \ref{psd-lf-vn} and \ref{neu_sch2N}, the $x$-axis represents the real part of $\zeta$ and the $y$-axis its
imaginary part.\\
We represent the domain $\{\zeta\in\C,|\lambda_\ell(\zeta)|\leq 1,\forall \ell\}$ delimited
by the curves $\{\zeta\in\C,\lambda_\ell(\zeta)=e^{i\theta},\theta\in[0,2\pi[\}$ for $0\leq \ell \leq s$.\\
On Fig. \ref{RK&AB} we plotted such domains for the first four Runge-Kutta schemes and the
first five Adams-Bashforth schemes.
Actually the curves correspond to all the values of $\zeta$, symbol of the operator $\delta t\, F$, for which there exists an eigenvalue
with modulus equal to one.
For order four and five Adams-Bashforth schemes, the stability domains only correspond to the semi-disks located on the left
of the $(Oy)$ axis.
The loops on the right of the $(Oy)$ axis do not correspond to any stable domain.

\ 

We discretize the differential equation (\ref{diffequ}) with respect to the space variables.
We assume that $u(t,x)=\sum_{k} u_k(t) \varphi_{\delta x,k}(x)\in V_{\delta x}$, where $\delta x$ is a
parameter corresponding to the space step. If $\overline u = (u_k)_k$, we obtain a discretized version
of (\ref{diffequ})
\begin{equation}
\label{diffequ_sv}
\partial_t \overline u = M_{\delta x}\,\overline u,\quad \overline u(0,\cdot)=\overline{u_0}
\end{equation}
with, for instance, $M_{\delta x}=\P_{\delta x}F$ where $\P_{\delta x}$ denotes the orthogonal projector
on $V_{\delta x}$.
To ensure the stability of the simulation, the stability domain with \emph{thick line}
must include the spectrum $Sp=\{\lambda(M_{\delta x})\}$ of the matrix $M_{\delta x}$.
The \emph{thick line} is due to the term $C\delta t$ in $\sup_{\xi\in Sp}\rho(M(\delta t \xi))\leq 1+C\delta t$
where $M$ from Eq. (\ref{defMsigma}) depends on the temporal scheme. 

\ 

The behavior of the stability domain along the $(Oy)$ axis indicates how the scheme will be stable under the
condition $\rho(M(\zeta))\leq 1+C\delta t$ --which gives more relevant stability conditions
than the more classical $\rho(M(\zeta))\leq 1$-- for convection-dominated problems.
The next parts of our study will be dedicated to finding precise stability conditions on $\delta t$ and $\delta x$
in the frame of von Neumann stability.

\section{Stability conditions for the transport equation}
\label{transp_eq}

The von Neumann stability analysis for the transport equation presents some subtleties which explain why Runge-Kutta
order two and Adams-Bashforth order two schemes are still used in numerical fluid dynamics although the transport
operator $ia\xi$ for fixed $a\in\R^*$ and $\xi\in[-\frac{\pi}{\delta x},\frac{\pi}{\delta x}]$,
is located outside the stability domains of these schemes.
In this section, we show that what matters is the behavior of the stability domain along the $(Oy)$ axis.

\subsection{Accurate theoretical stability condition}

Let us consider the most basic transport equation:
\begin{equation}
\label{transp}
\partial_t u+a\,\partial_x u=0,\qquad {\rm with}~~u:\R_+\times \R\to \R,\, (t,x)\mapsto u(t,x).
\end{equation}
Since $\widehat{f'}(\xi)=i\xi\widehat{f}(\xi)$, the symbol of the operator $F\,u=-a\,\partial_x u$ is equal to: $\sigma(\xi)=-i\,a\,\xi$.\\
As explained in the previous section, considering an explicit scheme (\ref{mixtype}), taking $u_n(x)=\phi_n\,e^{i\,\xi\,x}$, and setting $X_n$
as in equation (\ref{X_n}), we can write:
\begin{equation}
\label{ampliA}
X_{n+1}=A(\xi)X_{n}
\end{equation}
with $A$ a matrix whose coefficients are polynomials in $-i\,a\,\xi\,\delta t$.\\
In the case the numerical scheme is of Runge-Kutta type, then $A(\xi)$ is a polynomial:
\begin{equation}
\label{Aoxi}
A(\xi)=\sum_{\ell=0}^{s}\beta_\ell (-i\,a\,\xi)^\ell \delta t^{\ell}.
\end{equation}
The coefficients $(\beta_\ell)$ of this polynomial play an important role in our stability analysis.
In \cite{Pey00}, the polynomial $g(\zeta)=\sum_{\ell=0}^{s}\beta_\ell\zeta^\ell$ is called
the {\sl amplification factor}.
We are able to compute the norm of $A(\xi)$ explicitly:
\begin{equation}
\label{sumS_q}
|A(\xi)|^2=\sum_{\ell=0}^{s}S_\ell \delta t^{2\ell}a^{2\ell}\xi^{2\ell}
\end{equation}
with (assuming $\beta_{j}=0$ for $j>s$)
\begin{equation}
\label{calcSell}
S_\ell=\sum_{j=0}^{2\ell} (-1)^{\ell+j}\beta_j \beta_{2\ell-j}.
\end{equation}
The von Neumann stability condition $|A(\xi)|\leq 1+C\delta t$ for all frequencies $\xi$ remaining to the computational domain and for a given $C$,
implies that for $\xi\in[-\frac{1}{\delta x},\frac{1}{\delta x}]$, (usually the computational domain is rather $[-\frac{\pi}{\delta x},\frac{\pi}{\delta x}]$, but
we discard $\pi$ for simplicity):
\begin{equation}
\label{condSell}
\sum_{\ell=0}^{s}S_\ell \delta t^{2\ell}a^{2\ell}\xi^{2\ell}\leq 1+2C\delta t.
\end{equation}
For sake of consistency of the numerical scheme, $\beta_0=\beta_1=1$ so $S_0=\beta_0^2=1$. Then if $S_1=\dots=S_{r-1}=0$ and $S_r>0$ for
a given integer $r$, we can write for small $\delta t \xi$,
\begin{equation}
\label{Ared}
|A(\xi)|^2=1+S_r \delta t^{2r}a^{2r}\xi^{2r}+o(\delta t^{2r}\xi^{2r})
\end{equation}
with (\ref{condSell}) it implies $S_r \delta t^{2r}a^{2r}\delta x^{-2r}\leq 2C\delta t$,
so $\delta t^{2r}a^{2r}\xi^{2r}\to 0$ for $\delta t\to 0$ implies $\delta t\xi=o(1)$ i.e. as $\xi\sim 1/\delta x$,
we must have $\delta t=o(\delta x)$.
Hence the equation (\ref{Ared}) is valid for all the computational domain $[-\frac{1}{\delta x},\frac{1}{\delta x}]$.
And the stability condition (\ref{condSell}) is reduced to:
\begin{equation}
S_r \delta t^{2r}a^{2r}\delta x^{-2r}\leq 2C\delta t
\end{equation}
i.e.
\begin{equation}
\label{CFLdtdx}
\delta t\leq \left( \frac{2C}{S_r}\right)^{\frac{1}{2r-1}} \left( \frac{\delta x}{a} \right)^{\frac{2r}{2r-1}}.
\end{equation}

\ 

This surprising stability condition is directly linked to the tangency
of the stability domain $\mathcal{S}=\{\zeta\in\C,\max_i|\lambda_i(\zeta)|\leq 1\}$ to the vertical axis $(Oy)$.
Actually, we have the following theorem:
\begin{theorem}[Thick Line Stability Theorem]
\label{thtang}
Consider a numerical time integration of type (\ref{mixtype}) with stability domain $\mathcal{S}$
bounded near zero by the parameterized curve $\{\zeta(\theta),\theta\in\mathcal{V}(0)\}$ with $\mathcal{V}(0)$
a real neighborhood of zero.
If for some integer $r$, the Taylor expansion of $\zeta$ yields:
\begin{equation}
\zeta=i(\theta +o(\theta))+T_{2r}\theta^{2r}+o(\theta^{2r})
\end{equation}
with $T_{2r}<0$. Then the corresponding stability condition for the transport equation reads:
\begin{equation}
\label{eqFLST}
\delta t\leq \left( \frac{C}{-T_{2r}}\right)^{\frac{1}{2r-1}} \left( \frac{\delta x}{a} \right)^{\frac{2r}{2r-1}}.
\end{equation}
\end{theorem}
Remark that $T_{2r}=-\frac{S_r}{2}$ where the quantity $S_r$ defined by (\ref{calcSell}) provides (\ref{CFLdtdx}).

\ 

\noindent {\sl Proof:}
The amplification factor $g(\zeta)$ is given by $g(\zeta)=\lambda_{\max}(A(\xi))$ with $|\lambda_{\max}(A(\xi))|=\max|\lambda_i(\zeta)|$
and $\{\lambda_i(\zeta)\}_{0\leq i \leq s}$ the eigenvalues of the matrix $A(\xi)$ from equation (\ref{ampliA}).
There is a finite number of eigenvalues.
Due to the polynomial form of the elements of the matrix $A(\xi)$, the eigenvalues can be written as holomorphic functions:
$\lambda_{i}(\zeta)=\sum_{k\geq 0}\beta_{\ell}^{(i)}\zeta^\ell$.
Among these eigenvalues $\lambda_i(\zeta)$, we consider the one such that $|\lambda_{i_0}(\zeta)|\geq|\lambda_i(\zeta)|$ for all indices $i$
and complex number $\zeta$ in a neighborhood of $0$.
This corresponds to the largest sequence $(S_0,S_1,\dots,S_r,\dots)$ with $S_\ell$ defined by (\ref{calcSell}), for the usual order relation
on sequences.
Then $g(\zeta)=\lambda_{i_0}(\zeta)$ is a holomorphic function in a neighborhood of $0$:
\begin{equation}
\label{gtang}
g(\zeta)=\sum_{\ell\geq 0}\beta_{\ell}\zeta^\ell=1+\zeta+\beta_2 \zeta^2+\dots+\beta_s\zeta^s+\dots
\end{equation}
where, for consistency reasons, we consider $\beta_0=\beta_1=1$.\\
We already know that for $(S_\ell)$ given by (\ref{calcSell}) satisfying $S_1=\dots=S_{r-1}=0$ and $S_r>0$, the CFL stability
condition (\ref{CFLdtdx}) applies. We show that this same condition provides the tangency of the stability domain
$\mathcal{S}$ to the $(Oy)$ axis at zero.\\
Near $O$, let $\zeta=p+i\,q$,
\begin{equation}
\label{gtangpq}
g(\zeta)=1+p+iq+\beta_2 (p+iq)^2+\dots+\beta_s(p+iq)^s+\dots
\end{equation}
with $p$ and $q$ independent variables close to zero. Then,
\begin{equation}
|g(\zeta)|^2=(1+p+\beta_2 (p^2-q^2)+\dots)^2+(q+2\beta_2\,p\,q+\dots)^2.
\end{equation}
Looking for the first significant terms of this sum makes $1$ and $2p$ appear.
But we do not know which is the lowest power of $q$ existing in this sum.
Nevertheless, all the terms $p^\ell q^j$ with $\ell,j\geq 1$ and $p^\ell$ with $\ell\geq 2$ are negligible
with respect to $p$, so we have from (\ref{gtangpq})
\begin{eqnarray}
\nonumber |g(\zeta)|^2&=&|1+p+iq+\beta_2 (iq)^2+\dots+\beta_s(iq)^s+\dots|^2+o(p)\\
&=&1+2p+S_r\,q^{2r}+o(p)+o(q^{2r})
\end{eqnarray}
with $q^{2r}$ the lowest power of $q$ with nonzero coefficient $S_r$ (given by (\ref{calcSell})).\\
As a result, the curve $\{|g(\zeta)|=1\}$ is approximated by $p=-\frac{S_r}{2}q^{2r}$ near the origin.
This curve can also be parameterized by $\theta\in [-\pi,\pi]$ in $\{g(\zeta)=e^{i\theta},\theta\in [-\pi,\pi]\}$.
For multistep schemes (see Sec. \ref{secABsch}), it is convenient to express $\zeta$ as a function of $\theta$ and to
write it as a Taylor series:
\begin{equation}
\label{eqsumTell}
\zeta=\sum_{\ell\geq 1} i T_{2\ell-1}\theta^{2\ell-1}+T_{2\ell} \theta^{2\ell}.
\end{equation}
Then, for $\theta\in\R$ close to $0$
\begin{equation}
\label{expandzeta}
\zeta=i(\theta+o(\theta))-\frac{S_r}{2}\theta^{2r}+o(\theta^{2r}),
\end{equation}
so $T_{2r}=-\frac{S_r}{2}$. Hence this tangency implies the CFL (\ref{CFLdtdx}).
\begin{flushright}
$\Box$
\end{flushright}

\begin{theorem}
\label{theoO_CFL}
An order $2p$ numerical time integration applied to the transport equation is, at worst, stable under
the CFL-like condition:
\begin{equation}
\label{CFL2}
\delta t\leq C\, \left( \frac{\delta x}{a}\right)^{\frac{2p+2}{2p+1}}
\end{equation}
\end{theorem}
\emph{Proof:} 
For an order $2p$ scheme, we have:
\begin{equation}
\label{o2pschem}
u_{n+1}=u_n+\delta t\, \partial_t u_n+\frac{\delta t^2}{2}\partial_{t}^2 u_n+\dots
+\frac{\delta t^{2p}}{(2p)!}\partial_t^{2p} u_n+o(\delta t^{2p})
\end{equation}
The transport operator $F$ commutes with $\partial_t$.
So iterating $\partial_t u=F\,u$ we obtain $\partial_t^\ell u_n=F^\ell (u_n)$.
Hence equation (\ref{o2pschem}) yields the amplification factor:
\begin{equation}
\label{sum_taylor}
g(\zeta)=1+\zeta+\frac{\zeta^2}{2}+\dots+\frac{\zeta^{2p}}{(2p)!}+o(\zeta^{2p})
\end{equation}
with $o()$ gathering the negligible terms under the condition $\delta t=o(\delta x)$.
In this case, the $(\beta_\ell)$ of equation (\ref{gtang}) are given by $\beta_\ell=\frac{1}{\ell!}$.
Then for $q\in[1,p]$, the coefficients $S_q$ of the sum (\ref{sumS_q}) are given by:
\begin{equation}
S_q=\sum_{\ell=0}^{2q}(-1)^{(q-\ell)}\frac{1}{\ell!}~\frac{1}{(2q-\ell)!}=\frac{(-1)^q}{(2q)!}\sum_{\ell=0}^{2q}{\rm C}_{2q}^\ell(-1)^\ell=0
\end{equation}
Hence, in the worst case regarding the stability, the first nonzero significant term in the sum (\ref{sumS_q})
is $S_{p+1} \delta t^{2p+2}a^{2p+2}\xi^{2p+2}$ with $S_{p+1}>0$ implying the stability condition
(\ref{CFL2}). If $S_{p+1}<0$ then a linear CFL condition is sufficient.

\subsection{Examples with usual schemes}

We apply our analysis to some popular schemes in fluid dynamics.
This provides the following stability conditions for some of the most used schemes for transport
problem Eq. (\ref{transp}):
\begin{itemize}
\item The simplest example is the Euler explicit scheme, order one in time:
\begin{equation}
\label{EE}
u_{n+1}= u_{n} - \delta t ~a\partial_x u_n.
\end{equation}
For this scheme, $g(\zeta)=1+\zeta$ so $r=1$, $S_1=1$ and we find the stability condition:
\begin{equation}
\label{cflEE}
\delta t \leq 2 C \left(\frac{\delta x}{a} \right)^2.
\end{equation}
\item An improved version of this scheme allows us to construct an order two centered scheme:
\begin{equation}
\label{SC2}
\left\{
\begin{array}{ll}
u_{n+1/2}=& u_{n} - \frac{\delta t}{2} a \partial_x u_n\\
\vspace{-0.3cm} \\
u_{n+1}=& u_{n} - \delta t~ a\partial_x u_{n+1/2}
\end{array}
\right. .
\end{equation}
For this scheme, $g(\zeta)=1+\zeta+\frac{1}{2}\zeta^2$ so $r=2$ because $S_1=0$ and $S_2=\frac{1}{4}$.
Compared to the previous case, the stability is improved:
\begin{equation}
\label{cflSC2}
\delta t \leq 2 C^{1/3} \left(\frac{\delta x}{a} \right)^{4/3}.
\end{equation}
\item For Runge-Kutta scheme of order 4, we have:
\begin{equation}
\label{RK4}
\left\{
\begin{array}{ll}
u_{n(1)}=& u_{n} - \frac{\delta t}{2} a\partial_x u_n\\
\vspace{-0.3cm} \\
u_{n(2)}=& u_{n} - \frac{\delta t}{2} a\partial_x u_{n(1)}\\
\vspace{-0.3cm} \\
u_{n(3)}=& u_{n} - \delta t ~a\partial_x u_{n(2)}\\
\vspace{-0.3cm} \\
u_{n+1}=& u_{n} - \frac{\delta t}{6} a\partial_x u_{n}
- \frac{\delta t}{3} a\partial_x u_{n(1)}
- \frac{\delta t}{3} a\partial_x u_{n(2)}
  - \frac{\delta t}{6} a\partial_x u_{n(3)}
\end{array}
\right..
\end{equation}
From the amplification factor $g(\zeta)=1+\zeta+\frac{\zeta^2}{2}+\frac{\zeta^3}{6}+\frac{\zeta^4}{24}$ we infer
\begin{equation}
S_1=S_2=0\quad{\rm and}\quad S_3=-\frac{1}{72},\,S_4=\frac{1}{576}.
\end{equation}
As $S_3<0$, our study doesn't apply to this case, and the stability domain, Fig. \ref{RK&AB},
indicates that a classical linear CFL condition has to be satisfied.
\item The order 5 Runge-Kutta scheme from \cite{CM92} page 115 provides the
amplification factor $g(\zeta)=1+\zeta+\frac{\zeta^2}{2}+\frac{\zeta^3}{6}+\frac{\zeta^4}{24}+\frac{\zeta^5}{120}+\frac{\zeta^6}{1280}$.
Therefore it is stable under the condition:
\begin{equation}
\delta t\leq \left( \frac{11520}{7} \right)^{1/5} C^{1/5} \left( \frac{\delta x}{a} \right)^{6/5},~~~~~~\left( \frac{11520}{7} \right)^{1/5}\sim 4.398~.
\end{equation}
\item The order two Adams-Bashforth scheme goes as follows:
\begin{equation}
\label{AB2}
u_{n+1}= u_{n} - \frac{3}{2} \delta t ~a\partial_x u_{n}
+\frac{1}{2} \delta t ~a\partial_x u_{n-1}.
\end{equation}
So, according to Sec. \ref{VNstab}, we consider $X_n=\left[  \begin{array}{c} u_n \\ u_{n-1}\end{array}\right]$, and
we apply the numerical scheme to a pure Fourier mode $\phi_n\,e^{i\xi\, x}$. We obtain:
\begin{equation}
X_{n+1}=\left[  \begin{array}{cc} 1+\frac{3}{2}\zeta & -\frac{\zeta}{2} \\ 1 & 0 \end{array}\right] X_n
\end{equation}
with $\zeta=-ia\,\delta t\,\xi$.\\
We compute the eigenvalue of this $2\times 2$ matrix, the characteristic polynomial is given by
$\chi(Y)=Y^2-(1+\frac{3}{2}\zeta)Y+\frac{\zeta}{2}$. Owing to the fact that $\delta t=o(\delta x)$,
we have $\zeta \to 0$. An expansion of the larger eigenvalue $Y_0$ in terms of powers of $\zeta$ provides
\begin{equation}
Y_0=1+\zeta+\frac{\zeta^2}{2}-\frac{\zeta^3}{4}-\frac{\zeta^4}{8}+o(\zeta^4)
\end{equation}
With $\zeta=-ia\frac{\delta t}{\delta x}$, we obtain
\begin{equation}
|Y_0|=1+\frac{1}{4}a^4\frac{\delta t^4}{\delta x^4}+o(\frac{\delta t^4}{\delta x^4})
\end{equation}
As we want $|Y_0|\leq 1+C\delta t$, this drives to the following stability condition:
\begin{equation}
\label{cflAB2}
\delta t \leq 2^{2/3} C^{1/3} \left(\frac{\delta x}{a} \right)^{4/3}
\end{equation}
\end{itemize}

\ 

Therefore, two popular second order schemes, Runge Kutta two (RK2) and Adams-Bashforth two (AB2) require
a CFL-like condition: $\delta t \leq C \delta x^{4/3}$.
The $\delta t_{\rm max}$ is $2^{1/3}$ larger for RK2 than for AB2, but RK2 necessitates twice more computations than AB2
for each time step.
So, regarding only the stability, AB2 is $2^{2/3}$ cheaper than RK2.

\begin{figure}
\begin{center}
\centerline{\includegraphics[scale=0.5]{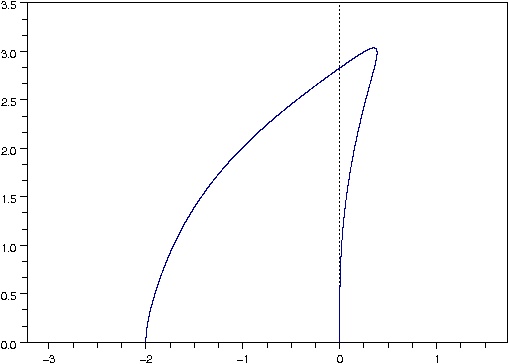}}
\caption{\label{psd-lf-vn} Von Neumann stability domain for the pseudo-Leap-Frog scheme equation (\ref{psd-lf}).}
\end{center}
\end{figure}

Not all the second order numerical schemes need to satisfy a $4/3$-CFL condition.
For instance, the Leap-Frog scheme calls a usual linear CFL stability condition.
The following second order scheme is also stable under a linear CFL condition:
\begin{equation}
\label{psd-lf}
u_{n+1}=u_n+\delta t \,a\partial_x\left(\frac{u_n+u_{n-1}}{2}+\delta t \,a\partial_x u_n\right)
\end{equation}
Its stability domain is drawn in Fig. \ref{psd-lf-vn}.
The fact that $r=2$ with $S_2<0$ in Eq. (\ref{expandzeta}) is reflected by a tangent to $(Oy)$ oriented to the right.

\subsection{Effect of the space discretization}
\label{space_discr}

\begin{figure}
\begin{center}
\centerline{\includegraphics[scale=0.45]{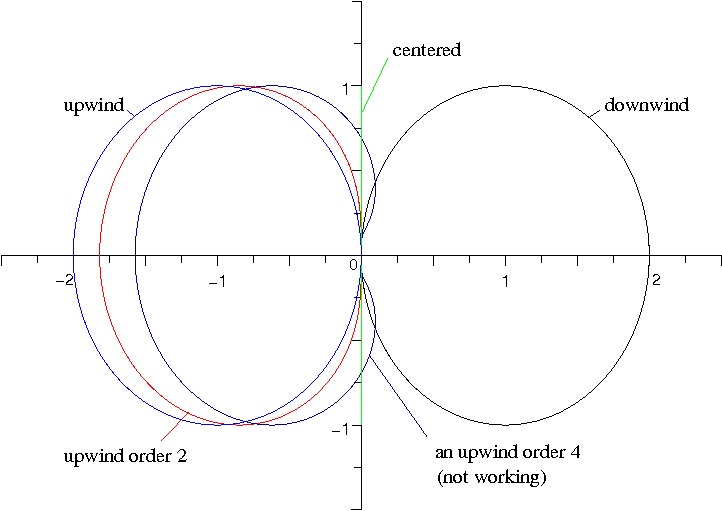}}
\caption{\label{stab_upcedo} Spectra of various finite difference schemes for space differentiation.
The numerical method is stable if the spectrum of the discretized differentiation fits into the domain of stability of
the temporal scheme (as those of Fig. \ref{RK&AB}). Here, the spectra have been normalized in order to have the same vertical size.}
\end{center}
\end{figure}

The space discretization impacts the stability condition (\ref{CFL2}) if it dissipates or creates energy, as do the upwind and downwind schemes.
Graphically this means that the spectra of these discretizations for the transport operator $F:u\mapsto -a\partial_x u$
are not contained in the $(Oy)$ axis, see Fig. \ref{stab_upcedo}.

In the frame of the von Neumann stability analysis we consider the function $u(x)=e^{i\xi x}$, for $\xi\in\R$.
Then, having a closer look at the three academic cases for finite differences, we obtain:
\begin{enumerate}
\item The downwind schemes are always unstable. For the first order downwind scheme (see Fig. \ref{stab_upcedo} for its spectrum)
\begin{equation}
a\frac{\partial u}{\partial x}+O(\delta x)=a\frac{u(x)-u(x-\delta x)}{\delta x}=a\,e^{i\xi x}\frac{-e^{-i\xi\delta x}+1}{\delta x}
\end{equation}
for $a<0$ provides the symbol $\sigma=-a\frac{-e^{-i\,\xi\,\delta x}+1}{\delta x}$ instead of $-i\,a\,\xi$ in formula (\ref{Aoxi}).
Combined with the Euler scheme for time integration, the amplification factor $G(\sigma)=1+\delta t\,\sigma$ becomes
\begin{equation}
|G|^2=1-2a\frac{\delta t}{\delta x}\left( 1+\frac{\delta t}{\delta x} \right)(1-\cos(\xi\,\delta x))
\end{equation}
and the error
\begin{equation}
\varepsilon_T\sim |G|^{\frac{T}{\delta t}}\varepsilon_0\sim e^{\frac{-2a}{\delta x}}\varepsilon_0
\end{equation}
goes unconditionally to $+\infty$.
\item The centered space discretizations satisfy $Sp\subset i\R$ where $Sp=\{\sigma(\xi),\xi\in[-\frac{\pi}{\delta x},\frac{\pi}{\delta x}]\}$.
They include most of the compact finite difference schemes.
For instance the usual centered scheme
\begin{equation}
\frac{\partial u}{\partial x}+O(\delta x^2)=\frac{u(x+\delta x)-u(x-\delta x)}{2\delta x}=e^{i\xi x}\frac{e^{i\xi\delta x}-e^{-i\xi\delta x}}{2\delta x}
\end{equation}
has a spectrum given by $\sigma=i\frac{\sin(\xi\,\delta x)}{\delta x}$ which goes along the $(Oy)$ axis, so its distance to the domain of
stability of the time scheme goes almost the same as in the spectral case. The stability results (\ref{eqFLST}) presented in this section apply fully to this case
with a constant $C$ which depends on the space discretization.
\item The upwind schemes can be unconditionally unstable if part of their spectra is located on the right side
of the $(Oy)$ axis (see the spectrum of the order four upwind scheme plotted on Fig. \ref{stab_upcedo}).
If their spectra remain in the left part of the complex plane, then the exponent in (\ref{eqFLST}) is modified in the following way:
assume that the domain of stability of the time discretization satisfies 
\begin{equation}
g(\theta)=i\theta+T_{2q}\theta^{2q}+o(i\theta)+o(\theta^{2q})
\end{equation}
in a neighborhood of $0$, with $T_{2q}<0$ and $q\in\N$;
assume that the spectrum of the discretized derivative satisfies
\begin{equation}
\sigma(\theta)=i\theta+V_{2p}\theta^{2p}+o(i\theta)+o(\theta^{2p})
\end{equation}
with $V_{2p}<0$ (i.e. upwind scheme) and $p\in\N$. Then the Thick Line Stability condition (\ref{eqFLST}) becomes:
\begin{itemize}
\item the CFL condition $\delta t\leq C\delta x$, with $C$ a constant independent of $\delta t$ and $\delta x$ if $p\leq q$,
\item a mitigation of the nonlinear condition (\ref{CFLdtdx})
\begin{equation}
\delta t\leq C\delta x^{\frac{q(2p-1)}{p(2q-1)}}
\end{equation}
with $C$ a constant independent of $\delta t$ and $\delta x$, if $p\geq q$.
\end{itemize}
The details of the proofs and the numerical tests for these assertions will be presented in a further article.
Remark that the case $p= +\infty$ (i.e. switching to a centered finite difference scheme) makes the condition (\ref{eqFLST}) appear.
\end{enumerate}

\section{Numerical experiment with the Burgers equation}
\label{numBurgers}
In order to test our assertions, we proceed to a numerical experiment with the inviscid Burgers equation.
Although this is a nonlinear equation, we choose initial conditions such that it assimilates to a transport equation: the sinusoidal part
represents only 1\% of the transport amplitude. So, technically regarding the stability, it behaves like
a transport equation.
Then the various Fourier modes are \emph{naturally} activated during the experiment.
\begin{equation}\label{burgers_eqn_per}
\partial_t u + u \partial_x u = 0\quad {\rm for}~~(t,x)\in[0,T]\times\T, \quad {\rm and}~~u(0,\cdot)=u_0.
\end{equation}
In Sec. \ref{CFLr} and \ref{secABsch} we show numerical evidence that stability conditions (\ref{cflEE}), (\ref{cflSC2}), (\ref{cflAB2}) and (\ref{CFL2})
hold for this problem (replacing $a$ by $\|u\|_{L^{\infty}}$).

We solve equation (\ref{burgers_eqn_per}) numerically using a Fourier pseudo-spectral method \cite{CHQZ88}.
The scheme is de-aliased by truncation.
Most of the time integration methods presented in this paper are tested on this classical basic problem.

The initial condition for the numerical experiment is
$u_0(x)=10-0.1\sin(\pi x)$, in a periodic domain $x\in\Omega=[-1,1]$.
For $t < t_{max} = 10/\pi$ the equation admits a smooth exact solution $u(x,t)=u_0(a)$, where $a=a(x,t)$
is solution of the equation $a-x+u_0(a)t=0$.
However, the numerical solution is only sought for $t \in [0,1]$, in order to satisfy some regularity
requirements on the solution (see proposition \ref{thRn}).
To determine the admissibility of the numerical solution, we apply a criterion based on the total variation
norm (which is expected to be constant):
the numerical solution $u_n$ has to satisfy $||u_n||_{TV}\leq K||u_0(x)||_{TV}$ with $K=1.1$ for all $n$ such that $n\delta t\leq T=1$.

\begin{figure}
\begin{center}
\includegraphics[scale=0.7]{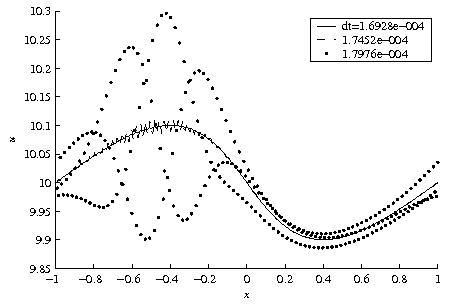}
\caption{\label{fig_explosion_RK2_N256} Numerical solution obtained at time $T=1$ for three
different time steps: $0.97\,\delta t_{\rm max}$, $\delta t_{\rm max}$ and $1.03\,\delta t_{\rm max}$ for $N=256$
with RK2 (order two Runge-Kutta) with a $\delta t_{\rm max}$ corresponding to $K=5$ i.e. $||u_T||_{TV}= 5||u_0(x)||_{TV}$.}
\end{center}
\end{figure}

The $\delta t_{\rm max}$ we compute, has very little dependence on the divergence criterion $K$.
Actually, below $\delta t_{\max}$ (97\%), the numerical solution shows no spurious oscillations,
while above it (103\%), these oscillations create some kind of explosion destroying the profile
of the solution completely, see Fig.~\ref{fig_explosion_RK2_N256}.

\begin{figure}
\centerline{
\includegraphics[width=13cm,clip]{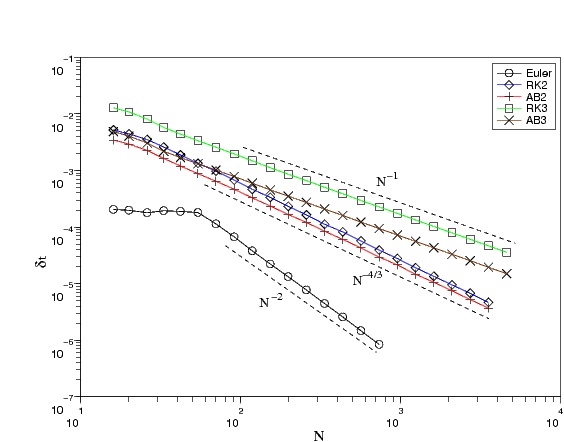}
}
\caption{\label{fig_stab_burg} Maximal time step $\delta t_{\max}$ depending on the number of points $N$ for
Runge--Kutta schemes and Adams--Bashforth schemes, obtained experimentally.}
\end{figure}

The computations are performed for different numbers of grid points, $16\leq N\leq 7758$.
For each $N$, we find $\delta t_{max}$ by dichotomy with a $0.5\%$ accuracy.
The results are represented as $\delta t_{max} (N)$ curves in Fig.~\ref{fig_stab_burg}.
They evidence the theoretically predicted power law $\delta t_{max}=CN^{\alpha}$
when the number of grid points is sufficiently large.
The explicit Euler scheme displays $\alpha=-2$ slope in Log-Log scale.
The two curves corresponding to the second-order schemes asymptotically both show an asymptotic
slope equal to $CN^{-\frac{4}{3}}$, but
the constant $C$ is $2^{1/3}$ times larger for the Runge--Kutta
scheme.

When the order is increasing to 3 and 4 for Runge--Kutta schemes and Adams--Bashforth schemes,
the slope equals $-1$.
But, while  the constant $C$ increases with the order for Runge--Kutta (yielding a larger stability domain),
it diminishes for Adams--Bashforth schemes with the increasing order
(see Fig. \ref{RK&AB} or \textit{e.g.}, \cite{CHQZ88}).

\section{Simple $\mathbf{2N}$-storage numerical schemes with ``shrinking CFL'' stability conditions}
\label{CFLr}

In order to illustrate the phenomenon presented in Sec. \ref{transp_eq}, we construct
numerical schemes having stability conditions of the type $\delta t \leq C \left(\frac{\delta x}{a}\right)^{\frac{2r}{2r-1}}$,
and which only necessitate two time levels to be stored in the computer memory.
Four of the five schemes presented here need to satisfy this stability condition with
exponents $2r/(2r-1)$ different from $1$: $2$, $\frac{4}{3}$, $\frac{6}{5}$ and $\frac{8}{7}$.
All of these numerical schemes are of order two, so they show relatively poor consistency given the number of
intermediate steps.
Other efficient low storage schemes can be found in \cite{Pey00} and \cite{GST01}.

To solve the equation
\begin{equation}
\partial_t u=F(u),
\end{equation}
let us consider the following family of schemes:
\begin{eqnarray}
\nonumber u_{(0)}&=&u_n\\
\nonumber u_{(1)}&=&u_{n}+\alpha_p \delta t F(u_{(0)})\\
\nonumber \dots & & \\
 u_{(\ell)}&=&u_{n}+\alpha_{p-\ell} \delta t F(u_{(\ell-1)})\\
\nonumber \dots & & \\
\nonumber u_{n+1}&=&u_{n}+\alpha_{1} \delta t F(u_{(p-1)})
\end{eqnarray}
These can also be written:
\begin{equation}
\label{simpl2Nsto}
u_{n+1}=u_n+\alpha_{1} \delta t F(u_n+\alpha_{2} \delta t F(u_n+\alpha_{3} \delta t F(u_n+\dots+\alpha_{p-1} \delta t F(u_n+\alpha_p\delta t F(u_n))\dots)))
\end{equation}
If $F$ is {\bf linear}, this corresponds to
\begin{equation}
\label{eq_sum_beta}
u_{n+1}=u_n+\beta_{1} \delta t Fu_n+\beta_{2} \delta t^2 F^2u_n+\beta_{3} \delta t^3 F^3u_n+\dots+\beta_{p} \delta t^p F^pu_n
\end{equation}
with $\beta_m=\prod_{\ell=1}^m \alpha_\ell$. Owing to $F^\ell u=\partial_t^{\ell} u$,
\begin{equation}
u_{n+1}=u_n+\beta_{1} \delta t \partial_tu_n+\beta_{2} \delta t^2 \partial_t^2u_n+\beta_{3} \delta t^3 \partial_t^3u_n+\dots+\beta_{p} \delta t^p \partial_t^pu_n.
\end{equation}
Here we recognize an expansion similar to the Taylor expansion of the function $u_n$, and we are able to tell exactly the
order of the scheme for linear equations by comparing the coefficients $\beta_\ell$ with those of the Taylor expansion which is provided by:
\begin{equation}
u_{n+1}=\sum_{\ell=0}^{+\infty} \frac{\delta t^\ell}{\ell!} \partial_t^\ell u_n
=u_n+\delta t \partial_tu_n+\frac{1}{2} \delta t^2 \partial_t^2u_n+\frac{1}{6} \delta t^3 \partial_t^3u_n+\dots
\end{equation}
and the smallest $\ell$ such that $\beta_{\ell+1}\neq 1/(\ell+1)!$ indicates the order of the scheme.
Remark that this holds only if $F$ is linear or if the order of the scheme is less or equal to two.
The interest of such schemes is that the coefficients $\alpha_\ell$ are easily deduced from the $\beta_\ell$.

\ 

We assume that $F$ is a convection operator.
Using the stability analysis Sec. \ref{transp_eq}, we know that the values of
\begin{equation}
S_\ell=\beta_\ell^2-2\beta_{\ell-1}\beta_{\ell+1}+2\beta_{\ell-2}\beta_{\ell+2}-\dots
\end{equation}
provide the stability condition.
We verify the validity of this stability condition using the numerical test from Sec. \ref{numBurgers}.

For a given $p$ in (\ref{simpl2Nsto}), maximizing the number of $S_\ell$ equal to zero leads to the following schemes
of order two --except the first one-- and stability conditions:
\begin{itemize}
\item with $\beta_1=1$ and $\beta_\ell=0$ for $\ell \geq 2$, this is the Euler explicit scheme:
\begin{equation}
u_{n+1}=u_n+\delta t Fu_n
\end{equation}
$\beta_2\neq \frac{1}{2}$ so it is of order $1$, and $S_1=1$ implies
$\delta t \leq 2C \left(\frac{\delta x}{a}\right)^2$.

\item with $\beta_1=1$, $\beta_2=1/2$ and $\beta_\ell=0$ for $\ell \geq 3$, this is a second order Runge-Kutta scheme:
\begin{equation}
u_{n+1}=u_n+\delta t F(u_n+\frac{1}{2} \delta t Fu_n)
\end{equation}
$\beta_3\neq \frac{1}{6}$ so it is of order $2$, and $S_2=1/4$ implies (\ref{cflSC2})
$\delta t \leq 2\,C^{1/3} \left(\frac{\delta x}{a}\right)^{4/3}$.

\item with $\beta_1=1$, $\beta_2=1/2$, $\beta_3=1/8$ and $\beta_\ell=0$ for $\ell \geq 4$, it is an order
two numerical scheme ($\beta_3\neq 1/6$),
\begin{equation}
{\rm (scheme~3)}\quad u_{n+1}=u_n+\delta t F(u_n+\frac{1}{2} \delta t F(u_n+\frac{1}{4} \delta t Fu_n))
\end{equation}
and as $S_1=S_2=0$ and $S_3=1/64$, we have the stability condition 
\begin{equation}
\label{cflschRK65}
\delta t \leq 2^{7/5}\,C^{1/5} \left(\frac{\delta x}{a}\right)^{6/5}.
\end{equation}

\item the schemes verifying $\beta_\ell=0$ for $\ell \geq 5$, and $S_1=S_2=S_3=0$ are given by
$\beta_1=1$, $\beta_2=1/2$, $\beta_3=\frac{2\pm \sqrt{2}}{4}$ and $\beta_4=\frac{3\pm 2\sqrt{2}}{8}$.
If we choose the minus sign for $\beta_3$ and $\beta_4$, this means:
\begin{equation}
{\rm (scheme~4)}\quad u_{n+1}=u_n+\delta t F(u_n+\frac{1}{2} \delta t F(u_n+\frac{2-\sqrt{2}}{2} \delta t F(u_n+\frac{2-\sqrt{2}}{4} \delta t Fu_n)))
\end{equation}
It is a second order scheme and has to satisfy the CFL-like stability condition
\begin{equation}
\label{cflschRK87}
\delta t \leq \left( \frac{2C}{\beta_4^2} \right)^{1/7} \left(\frac{\delta x}{a}\right)^{8/7}.
\end{equation}

\item in the general case, we consider $\beta_0=1,\beta_1=1,\beta_2,\dots,\beta_m$, and for $1\leq\ell\leq m-1$, $S_\ell=0$.
As $S_m=\beta_m^2>0$, the schemes resulting from this system of equations has to satisfy:
\begin{equation}
\label{withbeta}
\delta t \leq \left( \frac{2C}{\beta_m^2} \right)^{\frac{1}{2m-1}} \left(\frac{\delta x}{a}\right)^{\frac{2m}{2m-1}}.
\end{equation}
For $\beta_m$ positive and minimum, it results the constants indicated in Table \ref{table_beta}.
\begin{table}
\begin{center}
\begin{tabular}{|c|ccccccc|}
\hline
$m$ & ~~~~~~1~~~~~~ & ~~~~~~2~~~~~~ & ~~~~~~3~~~~~~ & ~~~~~~4~~~~~~ & ~~~~~~5~~~~~~ & ~~~~~~6~~~~~~ & ~~~~~~7~~~~~~ \\
\hline &&&&&&&\\
$\beta_m$ & 1 & $\frac{1}{2}$ & $\frac{1}{8}$ &  $\frac{3-2\sqrt{2}}{8}$ & $\frac{5\sqrt{5}-11}{64}$ & $\frac{26-15\sqrt{3}}{16}$ & $-\frac{1}{64}+\frac{7\alpha}{128}$  \\
&&&&&&&\\
\hline &&&&&&&\\
$\left(\frac{1}{\beta_m^2}\right)^\frac{1}{2m-1}$  & 1 & 1.587\dots & 2.297\dots & 2.997\dots & 3.687\dots & 3.395\dots & 5.045\dots    \\&&&&&&&\\
\hline \end{tabular}
\caption{\label{table_beta} Coefficients $\beta_m$ for different $m$. In the expression for $m=7$, $\alpha$ is a real solution of the equation $\alpha^3-9\alpha^2-\alpha+1=0$.}
\end{center}
\end{table}
\end{itemize}
On the other hand, if we impose the order to be 3 with five nonzero $\beta_\ell$, maximizing the number of $S_\ell$ equal to zero
provides $\beta_1=1$, $\beta_2=1/2$, $\beta_3=1/6$, $\beta_4=1/24$ and $\beta_5=1/144$. Hence it is written
\begin{equation}
{\rm (scheme~5)}\quad u_{n+1}=u_n+\delta t F(u_n+\frac{\delta t}{2} F(u_n+\frac{\delta t}{3} F(u_n+\frac{\delta t}{4} F(u_n+\frac{\delta t}{6} Fu_n))))
\end{equation}
As $S_4=\beta_4^2-2\beta_3\beta_5<0$, a classical linear CFL condition $\delta t\leq C\frac{\delta x}{a}$ applies.
Even, as $\beta_\ell=1/\ell!$ until $\ell=4$, this scheme is of order 4.

\begin{figure}
\centerline{\includegraphics[width=12cm,clip]{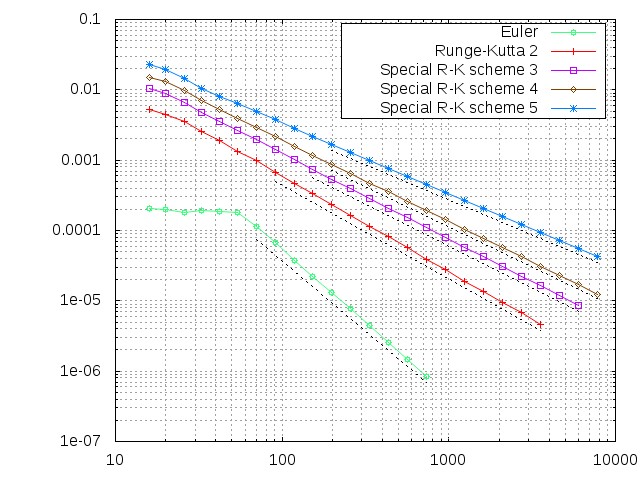}}
\caption{\label{stab_sch2N} Maximal time step ensuring stability for the test case of Sec. \ref{numBurgers} with the
Runge-Kutta schemes stable under the condition $\delta t \leq C \delta x^{\frac{2r}{2r-1}}$ for $r=1,2,3$ and $4$
whose slops we plot in parallel beneath the experimental results. $(Ox)$ axis
represents the number of points $N=\frac{2}{\delta x}$ and $(Oy)$ the maximal time step $\delta t_{\max}$ above which,
the numerical solution becomes unstable.}
\end{figure}

In Fig. \ref{stab_sch2N}, the slopes of stability condition on $\delta t$ issued from numerical experiments confirm our predictions for these schemes.

\begin{remark}
Maximizing the tangency of the stability domain to the $(Oy)$ axis is equivalent to optimizing
the energy conservation scale by scale. This explains why people simulating convection dominated
problems tend to prefer the Crank-Nicholson scheme
$u_{n+1}=u_n+\frac{\delta t}{2}(Fu_n+Fu_{n+1})$ (see \cite{CH08} for instance) whose stability
domain boundary coincides with the $(Oy)$ axis.
\end{remark}

\section{Adams-Bashforth schemes with ``shrinking CFL'' stability conditions}
\label{secABsch}
Let us consider an Adams-Bashforth scheme with coefficients $(\alpha_k)$:
\begin{equation}
\label{eqABsch}
u_{n+1}=u_{n}+\sum_{k=0}^{K}\alpha_k\,\delta t \,F(u_{n-k}).
\end{equation}
The order of scheme (\ref{eqABsch}) depends on the sums:
\begin{equation}
\label{sumordAB}
\Upsilon_\ell=\sum_{k=0}^{K} k^\ell\alpha_k,\qquad {\rm for}~0\leq \ell \leq K~. 
\end{equation}
The scheme has order $m$, iff for $0\leq \ell\leq m-1$, $\Upsilon_\ell=\frac{(-1)^\ell}{\ell+1}$ (see \cite{Pey00}).
Solving the system for $m=K+1$ provides the Adams-Bashforth scheme of order $K+1$ properly speaking.

The von Neumann stability domain is computed as indicated in Sec. \ref{VNstab}. Let
\begin{equation}
X_n=\left( \begin{array}{c}
u_n \\ u_{n-1} \\ \vdots \\ u_{n-K}
\end{array} \right)  \qquad {\rm and} \qquad \widehat{\delta t\,F(u_\ell)}=\zeta \widehat{u_\ell}~~{\rm with}~~\zeta\in\C.
\end{equation}
Then $\widehat{X_{n+1}}=M(\zeta)\widehat{X_n}$, with the matrix $M(\zeta)\in\mathcal{M}_{K+1}(\C)$ given by
\begin{equation}
\label{eqAbmat}
M(\zeta)=\left[  \begin{array}{ccccc} 1+\alpha_0\zeta & \alpha_1\zeta & \dots & \dots & \alpha_K\zeta \\
1 & 0 & \dots & \dots & 0 \\ 0 & \ddots & \ddots & & \vdots \\ \vdots & \ddots & \ddots & \ddots &\vdots
\\ 0 & \dots & 0 & 1 & 0 \end{array}\right].
\end{equation}
The characteristic polynomial is given by
\begin{equation}
\label{eqpolcar}
P(X)={\rm det}(X\,Id-M(\zeta))=X^{K+1}-X^K-\sum_{k=0}^{K}\alpha_k\zeta X^{K-k}.
\end{equation}
As the eigenvalues of this polynomial provide the multiplication factor of the scheme (\ref{eqABsch}),
the boundaries of the stability domain are therefore obtained by considering the curve
$\{\zeta\in\C~{\rm s.t.}~\exists\theta\in[-\pi,\pi],P(e^{i\theta})=0\}$, i.e.
\begin{equation}
\label{zetAB}
\zeta=\frac{e^{i\theta}-1}{\sum_{k=0}^{K}\alpha_k\,e^{-ik\theta}}, \qquad \theta\in[-\pi,\pi].
\end{equation}
According to Theorem \ref{thtang}, the stability condition depends on the tangency to the imaginary axis $(Oy)$
obtained for $\theta$ close to $0$.
Assuming order one at least (i.e. $\Upsilon_0=1$), a Taylor expansion of expression (\ref{zetAB}) provides:
\begin{equation}
\label{eqTayzet}
\zeta=\sum_{r\geq 1} i^r \theta^r \sum_{\begin{array}{c}
                                               p+q=r \\ p\geq 1,\,q\geq 0
                                              \end{array}}
\frac{(-1)^q}{p!}\sum_{\sum_{n\geq 1}n\kappa_n=q}(-1)^{\sum_{n\geq 1}\kappa_n} \frac{\left(\sum_{n\geq 1}\kappa_n\right)!}{\prod_{n\geq 1}\kappa_n!}
\prod_{n\geq 1}\left( \frac{\Upsilon_n}{n!} \right)^{\kappa_n}.
\end{equation}
The first two elements of this sum are given by
\begin{equation}
T_2=-\Upsilon_1-\frac{1}{2},\quad
T_4=\frac{1}{6} \Upsilon_1-\frac{1}{4}\Upsilon_2-\Upsilon_1\Upsilon_2+\frac{1}{6}\Upsilon_3+\frac{1}{2}\Upsilon_1^2+\Upsilon_1^3
\end{equation}
with $\Upsilon_\ell$ from (\ref{sumordAB}).

For a given $K$, maximizing the tangency to $(Oy)$ (i.e. the number $m$ s.t. $T_{2\ell}$ is equal to $0$ for $l<m$)
provides the following numerical schemes:
\begin{itemize}
 \item for $K=1$, $T_2=0$ implies $\alpha_0=\frac{3}{2}$ and $\alpha_1=-\frac{1}{2}$, i.e. Adams-Bashforth scheme
of order two.
As $T_4=-\frac{1}{4}$, it is stable under the condition (\ref{cflAB2}) $\delta t \leq 2^{2/3} C^{1/3} \left(\frac{\delta x}{a} \right)^{4/3}$.
 \item with three time steps, $T_2=T_4=0$ leads to the scheme we call (ABsch3) with
$\alpha_0=\frac{5}{3}$, $\alpha_1=-\frac{5}{6}$ and $\alpha_2=\frac{1}{6}$.
Given that $\Upsilon_1=-1/2$ and $\Upsilon_2=-1/6\neq 1/3$, it is of order two. And $T_6=-1/12$ induces the CFL condition
\begin{equation}
\label{cdAB65}
\delta t \leq 12^{1/5} C^{1/5} \left(\frac{\delta x}{a} \right)^{6/5}.
\end{equation}
 \item with four time steps, enforcing $T_2=T_4=T_6=0$ yields the scheme (ABsch4) with
\begin{equation}
(\alpha_0,\alpha_1,\alpha_2,\alpha_3)=(\frac{7}{4},-\frac{21}{20},\frac{7}{20},-\frac{1}{20})
\end{equation}
As $\Upsilon_1=-1/2$ and $\Upsilon_2=1/10\neq 1/3$, this is also a second order scheme. On the other hand, we have $T_8=-\frac{1}{40}$,
so this scheme is stable under the condition
\begin{equation}
\label{cdAB87}
\delta t \leq 40^{1/7} C^{1/7} \left(\frac{\delta x}{a} \right)^{8/7}.
\end{equation}
\end{itemize}

\begin{figure}
\centerline{\includegraphics[width=11cm,clip]{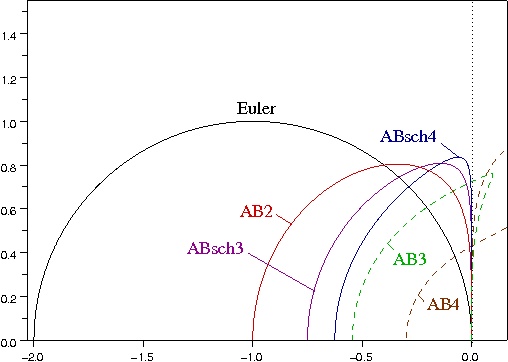}}
\caption{\label{neu_ABm} Von Neumann stability domains for modified Adams-Bashforth schemes maximizing the
tangency to the $(Oy)$ axis.}
\end{figure}

\begin{figure}
\centerline{\includegraphics[width=12cm,clip]{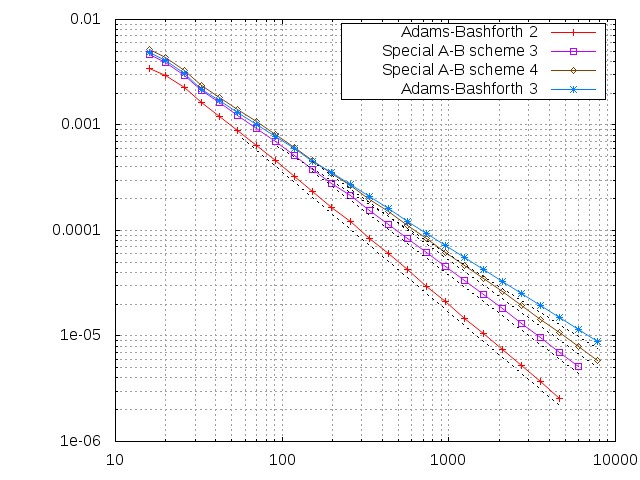}}
\caption{\label{stab_ABm} Maximal time step insuring the stability, obtained experimentally with the 1-D Burgers equation (\ref{burgers_eqn}) for
Adams-Bashforth schemes. It evidences three slopes: $\delta t_{max} = C \delta x^{2r/(2r-1)}$ with $r=2,3$ and $4$, plus a CFL condition slope,
plotted in parallel with dotted lines.
$(Ox)$ axis represents the number of points $N=\frac{2}{\delta x}$ and $(Oy)$ the maximal time step $\delta t_{\max}$ above which,
the numerical solution becomes unstable.}
\end{figure}

We plot the stability domains corresponding to these schemes on Fig. \ref{neu_ABm}, and verify our
stability predictions with the Burgers equation test Sec. \ref{numBurgers}. The results of these experiments
on Fig. \ref{stab_ABm} confirm the predicted stability conditions (\ref{cflAB2}), (\ref{cdAB65}) and (\ref{cdAB87}),
but less accurately than for Runge-Kutta schemes (\ref{cflSC2}), (\ref{cflschRK65}) and (\ref{cflschRK87}).

\section{Extension to some nonlinear equations}
\label{nlin_eq}
We show that these results extend to regular solutions to nonlinear
problems such as the incompressible Euler equations on a domain $\Omega$
bounded with walls, and scalar conservation laws.
We proceed in three steps with gradually increasing complexity:
\begin{itemize}
\item First we consider the transport equation with non-constant velocity on bounded domains.
Hence we step outside the strict frame of von Neumann stability analysis.
\item Then we study the simplest nonlinear equation involving transport: the 1D Burgers equation,
and we show that the previous results still hold true under a smoothness condition.
\item Then we transpose our results to the scalar conservation laws and the incompressible Euler
equations on a domain $\Omega$ possibly bounded by walls.
\end{itemize}

\subsection{Transport by a variable velocity}
\label{transpU}

The transport of a scalar $\theta$ by a divergence-free velocity $\u$ on an open set $\Omega\subset \R^d$
with regular boundaries satisfies the equation:
\begin{eqnarray}
\label{eqtranspU}
\partial_t \theta +\u(\x)\cdot \nabla\theta=0 ~{\rm for}~ \x\in\Omega,~t\in [0,T],\\
\nonumber{\rm with}\quad{\rm div}(\u)=0~{\rm for}~ \x\in\Omega, \quad \u(\x)\cdot {\bf n}=0\quad {\rm for}~~\x\in\partial\Omega.
\end{eqnarray}
In order to generalize the stability analysis to this case, we need the following lemma which corresponds
to $\v=(\theta,\dots,\theta)$ and $\w=(\varphi,\dots,\varphi)$ in lemma \ref{lemme}:
\begin{lemma}
\label{lemmeth}
Let $\theta,\varphi:\Omega\to \R$, $\u:\Omega\to \R^d$ such that ${\rm div}(\u)=0$ on $\Omega$, and $\u\cdot{\bf n}=0$
on $\partial \Omega$ then:
\begin{equation}
\label{orthtrans}
\langle \theta,\u\cdot\nabla\varphi\rangle_{L^2(\Omega)}=-\langle \u\cdot\nabla\theta,\varphi\rangle_{L^2(\Omega)}
\end{equation}
Equivalently, we have:
\begin{equation}
\label{orthtrans2}
\langle \theta,\u\cdot\nabla\theta\rangle_{L^2(\Omega)}=0
\end{equation}
\end{lemma}
The computations using the skew-symmetry relationship leads to the same stability conditions
as those relying on complex numbers in Sec. \ref{transp_eq} under the following assumptions
regarding the space discretization:

\begin{assumption}
\label{assu1}
The space discretization conserves the skew-symmetry of the equation, i.e. with the notations of equation (\ref{diffequ_sv})
\begin{equation}
\langle M_{\delta x}\overline \theta,\overline \theta \rangle=0,\qquad \textrm{which is equivalent to}
\quad \langle M_{\delta x}\overline \theta,\overline \varphi \rangle
= -\langle \overline \theta,M_{\delta x}\overline \varphi \rangle,\quad \forall \overline \theta,\overline \varphi\in V_{\delta x}.
\end{equation}
\end{assumption}
Such \emph{conservative} discretizations are presented in some computational fluid mechanics publications such as \cite{VV03} for instance.

\begin{assumption}
\label{assu2}
The discretization is sufficiently regular to enforce
\begin{equation}
\|M_{\delta x}\overline \theta\|\leq C \frac{\|\overline \theta\|}{\delta x},\quad \forall \overline \theta\in V_{\delta x}.
\end{equation}
\end{assumption}
This assumption is satisfied with $C\sim \|\u\|_{L^{\infty}}$ for almost all the discretizations.
One would need special properties to avoid this to happen.

Let $F(\theta)=\P_{\delta x}(\u\cdot\nabla\theta)$ with $\P_{\delta x}$ the orthogonal projector onto the space of discretization $V_{\delta x}$.
In the next sections Sec. \ref{burger}, \ref{cons-law} and \ref{euleri}, even if the operator $F$ is not linear,
we omit the projector $\P_{\delta x}$ since it does not change the computations we present
because it can be set or removed when needed:
\begin{equation}
\label{Pdx_inutile}
\langle \P_{\delta x} \theta,\P_{\delta x} \varphi \rangle=\langle \P_{\delta x} \theta,\varphi \rangle=\langle \theta,\P_{\delta x} \varphi \rangle,
\quad {\rm for}~\theta,\varphi\in V_{\delta x}.
\end{equation}
For the Runge-Kutta scheme (\ref{RKtype}), we find the following expression for $\theta_{n+1}$:
\begin{equation}
\label{th_n+1_Fi}
\theta_{n+1}=\sum_{i=0}^{k} \beta_{i}\delta t^iF^i(\theta_n)
\end{equation}
Starting from this expression and according to lemma \ref{lemmeth} along with assumption \ref{assu1},
\begin{equation}
\label{_FF}
\langle F^i(\theta_n),F^j(\theta_n)\rangle_{L^2(\Omega)}=\left\{
\begin{array}{ll} 0 &~~~~{\rm if}~~~~i+j=2\ell+1~~~~{\rm for}~~~~\ell\in\N \\
(-1)^{\ell-i}\|F^\ell(\theta_n)\|_{L^2}^2 &~~~~{\rm if}~~~~i+j=2\ell~~~~{\rm for}~~~~\ell\in\N
\end{array}\right..
\end{equation}
We compute the $L^2$ norm of $\theta_{n+1}$ as a function of the $L^2$ norm of $\theta_{n}$.
From (\ref{th_n+1_Fi}), (\ref{_FF}) and under the assumption \ref{assu1}, we have:
\begin{equation}
\label{sumS}
\|\theta_{n+1}\|_{L^2}^2=\sum_{\ell=0}^{k} S_{\ell}\,\delta t^{2\ell}\|F^\ell(\theta_n)\|_{L^2}^2
\end{equation}
with $(S_\ell)$ given by (\ref{calcSell}), i.e.
\begin{equation}
\label{Sell}
S_{\ell}=\sum_{j=-{\rm min}(\ell,k-\ell)}^{{\rm min}(\ell,k-\ell)}(-1)^j\beta_{\ell-j}\beta_{\ell+j}
\end{equation}
For consistency needs of the numerical scheme, we must have $S_0=1$. On the other hand let us suppose that
$S_1=S_2=\dots=S_{r-1}=0$ and $S_{r}>0$. Under the assumption \ref{assu2}, for $\theta_n\in V_{\delta x}$,
\begin{equation}
\|F^r(\theta_n)\|_{L^2}\leq \|{\bf u}\|_{L^\infty}^{r} \frac{\|\theta_n\|_{L^2}}{\delta x^r},
\end{equation}
and knowing that for $x\geq -1$,
\begin{equation}
\sqrt{1+x}\leq 1+\frac{x}{2},
\end{equation}
we derive from (\ref{sumS}):
{\small \begin{equation}
\label{expS}
\|\theta_{n+1}\|_{L^2}
\leq\left( 1+\frac{\delta t^{2r}}{\delta x^{2r}}S_r\|{\bf u}\|_{L^\infty}^{2r}+o(\delta t) \right)^{1/2} \|\theta_{n}\|_{L^2}
\leq\left( 1+\left(\frac{\delta t^{2r-1}S_r}{2\delta x^{2r}}\|{\bf u}\|_{L^\infty}^{2r}+o(1)\right)\delta t \right) \|\theta_{n}\|_{L^2}
\end{equation}}
where $o()$ gathers all the negligible terms.\\
Let us note $a=\|{\bf u}\|_{L^\infty}$, then the numerical scheme (\ref{RKtype}) is stable for small perturbations under the condition:
\begin{equation}
\label{CFL1}
\delta t\leq C \left(\frac{\delta x}{a}\right)^{\frac{2r}{2r-1}}.
\end{equation}
Hence the results obtained in the von Neumann stability framework remain valid in the case
of the convection by a variable velocity on a bounded domain. This is still a linear equation but outside
the von Neumann stability analysis framework which assumes a periodic or unbounded domain.

\subsection{The Burgers equation}
\label{burger}

In order to clarify the role of the nonlinearity and validate our analysis under smoothness
conditions on the solution, we have a look at the simplest nonlinear case, the one-dimensional
inviscid Burgers equation:
\begin{equation}\label{burgers_eqn}
\partial_t u + u \partial_x u = 0\quad {\rm for}~~(t,x)\in [0,T]\times\R,\quad u(0,\cdot)=u_0.
\end{equation}
In order to infer the numerical stability for this problem, we linearize it.
Assume $u_n$ is a discretized version of the solution $u$ in time and in space.
As proposed in \cite{GR96}, we consider a perturbed solution $u_n+\varepsilon_n$.
Under regularity assumptions on $u$, each time discretization will involve
a specific evolution equation on $\varepsilon_n$.

Actually, the small error $\varepsilon_n$ that we introduce corresponds to oscillations at the smallest scales
in space $V_{\delta x}$. This instability propagates and may increase at each time step.
In the following, we demonstrate that under CFL-like conditions similar to those of Sec. \ref{transp_eq}, the $L^2$ norm of
the small error $\varepsilon_n$ is amplified in a limited way:
\begin{equation}
\label{stabCdt}
\|\varepsilon_{n+1}\|_{L^2}\leq (1+C\delta t) \|\varepsilon_{n}\|_{L^2}
\end{equation}
where $C$ is a constant that neither depends on $\delta x$ nor on $\delta t$.
Thus, after an elapsed time $T$, the error increases at most exponentially as a function of the time:
\begin{eqnarray}
\|\varepsilon_{t_0+T}\|_{L^2}\leq
\left(1+C\delta t\right)^{T/\delta t} \|\varepsilon_{t_0}\|_{L^2}
\leq e^{CT} \|\varepsilon_{t_0}\|_{L^2}
\end{eqnarray}
As $\partial_t u =- u \partial_x u$,
$\partial_t^\ell u =\sum_{\alpha}\lambda_\alpha u^{\alpha_1}(\partial_xu)^{\alpha_2}\dots(\partial_x^{\ell-1}u)^{\alpha_{\ell-1}}
+(-1)^\ell u^\ell \partial_x^\ell u$,
so we remark a kind of equivalence between the space regularity and the time regularity.
If $\partial_x^\ell u\in L^\infty$ then $\partial_t^\ell u\in L^\infty$.
In the general case, for Runge-Kutta schemes (\ref{RKtype}), we have for $0\leq \ell\leq s$,
\begin{equation}
\label{u+e}
u_{(\ell)}+\varepsilon_{(\ell)}=\sum_{i=0}^{\ell-1}a_{\ell i}\,(u_{(i)}+\varepsilon_{(i)})
+\delta t \,\sum_{i=0}^{\ell-1}b_{\ell i}\,F(u_{(i)}+\varepsilon_{(i)})
\end{equation}
and $u_{n+1}+\varepsilon_{n+1}=u_{(s)}+\varepsilon_{(s)}$, so
\begin{equation}
\label{e_(l)}
\varepsilon_{(\ell)}=\sum_{i=0}^{\ell-1}a_{\ell i}\,\varepsilon_{(i)}
+\delta t \,\sum_{i=0}^{\ell-1}b_{\ell i}\,(F(u_{(i)}+\varepsilon_{(i)})-F(u_{(i)}))
\end{equation}
and $\varepsilon_{n+1}=\varepsilon_{(s)}$.

\begin{proposition}
\label{thRn}
Consider a solution $u$ of the Burgers equation (\ref{burgers_eqn}) $s$-times differentiable such that $\|\partial_x^s u\|_{L^{\infty}(\R\times[0,T])}<+\infty$.
Under the condition $\delta t=o(\delta x)$, a stability error $\varepsilon_{n+1}$ in the explicit scheme (\ref{e_(l)}), small enough at the initial
time: $\|\varepsilon_0\|_{L^{2}}=o(\delta x^{3/2})$ can be expressed as
\begin{equation}
\label{varepsL2}
\varepsilon_{n+1}
=\varepsilon_{n}+\sum_{i=1}^s \beta_{i}\, \delta t^i\, u_n^i\,\partial_x^i\varepsilon_{n}+\delta t\, \varepsilon_n\,\partial_x u_n+R_{n}
\end{equation}
with $\|R_{n}\|_{L^2}=o(\delta t \|\varepsilon_n\|_{L^2})$.
The coefficients $(\beta_i)$ derive from scheme (\ref{RKtype}) similarly as those in (\ref{Aoxi}).
\end{proposition}
{\sl proof:}
All the terms we have to deal with are projections in the space discretization $V_{\delta x}$.
In order to simplify the notation, we omit this projection that we assume orthogonal, as for the Galerkin methods \cite{ZS04,BEFZ88}.\\
We prove that $\varepsilon_{(\ell)}$ can be put under the form (\ref{varepsL2}) by recurrence on $\ell=0\dots s$.\\
As $\varepsilon_{(0)}=\varepsilon_n$, the assertion is true for $\ell=0$.\\
Let us assume the assertion true for $i$ from $0$ to $\ell-1$:
\begin{equation}
\label{varepsL2i}
\varepsilon_{(i)}
=\varepsilon_{n}+\sum_{j=1}^i \beta_{(i)j}\, \delta t^j\, u_n^j\,\partial_x^j\varepsilon_{n}+\alpha_{(i)}\,\delta t\, \varepsilon_n\,\partial_x u_n+R_{(i)}
\end{equation}
with $\|R_{(i)}\|_{L^2}=o(\delta t \|\varepsilon_n\|_{L^2})$.
The coefficients $(\beta_{(i)j})$ correspond to the partial step $i$ of the Runge-Kutta scheme distant by $\alpha_{(i)}\delta t$
from the time $n\delta t$. Remark that $\alpha_{(s)}=1$.
Then, given that $\sum_{i=0}^{\ell-1}a_{\ell i}=1$,
\begin{eqnarray}
\nonumber\varepsilon_{(\ell)}&=&\sum_{i=0}^{\ell-1}a_{\ell i}\,\varepsilon_{(i)}
+\delta t \,\sum_{i=0}^{\ell-1}b_{\ell i}\,(F(u_{(i)}+\varepsilon_{(i)})-F(u_{(i)})) \\
\label{calc_epsell} &=&\varepsilon_n +\sum_{i=0}^{\ell-1}a_{\ell i}\left(
\sum_{j=1}^i \beta_{(i)j} \delta t^j u_n^j\partial_x^j\varepsilon_{n}+\alpha_{(i)}\delta t \varepsilon_n\partial_x u_n+R_{(i)}
\right)\\
\nonumber &&+\delta t \,\sum_{i=0}^{\ell-1}b_{\ell i}\,\left(u_{(i)}\partial_x\varepsilon_{(i)}
+\varepsilon_{(i)}\partial_x u_{(i)}+\varepsilon_{(i)}\partial_x\varepsilon_{(i)}\right).
\end{eqnarray}
Knowing that
\begin{equation}
\partial_x \varepsilon_{(i)}=\partial_x\varepsilon_{n}+\sum_{j=1}^i \beta_{(i)j} \delta t^j \left( \partial_x(u_n^j)\partial_x^j\varepsilon_{n}
+u_n^j\partial_x^{j+1}\varepsilon_{n} \right)
+\alpha_{(i)}\delta t \left( \partial_x \varepsilon_{n}\partial_x u_{n}+ \varepsilon_{n}\partial_x^2 u_{n} \right)+\partial_x R_{(i)},
\end{equation}
then
{\small \begin{eqnarray}
\nonumber R_{(\ell)}&=&\sum_{i=0}^{\ell-1}a_{\ell i}\,R_{(i)}
+\delta t \,\sum_{i=0}^{\ell-1}b_{\ell i}\,\left(u_{(i)}
\left( \sum_{j=1}^i \beta_{(i)j} \delta t^j \partial_x(u_n^j)\partial_x^j\varepsilon_{n}
+\delta t \left( \partial_x \varepsilon_{n}\partial_x u_{n}+ \varepsilon_{n}\partial_x^2 u_{n} \right)+\partial_x R_{(\ell)}\right)\right. \\
\label{R_ell} &&+\varepsilon_{(i)}\partial_x\varepsilon_{(i)}+(\varepsilon_{(i)}\partial_x u_{(i)}-\varepsilon_{n}\partial_x u_{n})
+(u_{(i)}-u_n)\sum_{j=1}^i \beta_{(i)j} \delta t^j u_n^j\partial_x^{j+1}\varepsilon_{n}
\Bigg).
\end{eqnarray}}
Now, we need to show that these terms are $o(\delta t \|\varepsilon_{n}\|_{L^2})$.
According to assumption (\ref{varepsL2i})and due to $\delta t=o(\delta x)$, assumption \ref{assu2} provides
\begin{equation}
\|\delta t ^j\partial_x^j\varepsilon_n\|\leq \left( \frac{\delta t}{\delta x} \right)^j\|\varepsilon_{n}\|
\end{equation}
Hence $\varepsilon_{(i)}=(1+o(1))\varepsilon_{n}$
in the sense $\varepsilon_{(i)}=\varepsilon_{n}+\eta_{(i)}$ with $\|\eta_{(i)}\|_{L^2}=o(\|\varepsilon_{n}\|_{L^2})$
(the stability condition being $\|\varepsilon_{(i)}\|_{L^2}=(1+O(\delta t))\|\varepsilon_{n}\|_{L^2}$).
As we assumed $\|\varepsilon_n\|_{L^{2}}=o(\delta x^{3/2})$, then
\begin{equation}
\|\varepsilon_n\|_{L^{\infty}}\leq\frac{\|\varepsilon_n\|_{L^{2}}}{\delta x^{1/2}}=o(\delta x).
\end{equation}
As a result, the cross term $\delta t\varepsilon_{(i)}\partial_x\varepsilon_{(i)}$ satisfies
\begin{equation}
\|\delta t\varepsilon_{(i)}\partial_x\varepsilon_{(i)}\|_{L^2}\leq \frac{\delta t}{\delta x}\|\varepsilon_{(i)}\|_{L^{\infty}}\|\varepsilon_{(i)}\|_{L^{2}}
=o(\delta t\|\varepsilon_{n}\|_{L^{2}}).
\end{equation}
As for $i\leq s-1$,
\begin{equation}
u_{(i)}=u_n+\delta t B_{(i)}(u_n,\partial_x u_n,\dots,\partial_x^i u_n,\delta t)
\end{equation}
with $B$ a polynomial, $\|B\|_{L^\infty}$ is bounded, as well as $\|\partial_x B\|_{L^\infty}$ so
$\|u_{(i)}-u_n\|_{L^\infty}$=o(1) and $\|\partial_x u_{(i)}-\partial_x u_n\|_{L^\infty}$=o(1).
It allows us to replace $u_{(i)}$ by $u_n$ in the expansion (\ref{calc_epsell}), the difference going into
$R_{(\ell)}$, see (\ref{R_ell}).\\
Hence, using the fact that $\varepsilon_{(i)},R_{(i)}\in V_{\delta x}$ the discretization space,
$\|\partial_x^j \varepsilon_{(i)}\|_{L^2}\leq\frac{\|\varepsilon_{(i)}\|_{L^2}}{\delta x^j}$ and the same for $R_{(i)}$.
Let $r$ be an element of the sum $R_{(\ell)}$, then it satisfies:
\begin{equation}
\|r\|_{L^2}\leq \frac{\delta t^p}{\delta x^q}\tau(\|u_n\|_{L^\infty},\|\partial_x u_n\|_{L^\infty},\dots,\|\partial_x^\ell u_n\|_{L^\infty})\|\varepsilon_{n}\|_{L^2}
\end{equation}
with $\tau$ a polynomial, and $p\geq q+1$.\\
Given the fact that $\delta t=o(\delta x)$, we obtain that $\|R_{(\ell)}\|_{L^2}=o(\delta t\|\varepsilon_{n}\|_{L^2})$.
Using the recurrence, we obtain the result for $\ell=s$ i.e. for $\varepsilon_{n+1}$.

Actually, taking into account the orthogonality of $\varepsilon_n\partial_x\varepsilon_n$ with $\varepsilon_n$,
we can relax one of the assumptions i.e. it is sufficient to have $\|\varepsilon_0\|_{L^{2}}=o(\delta x)$,
and with the cancellations, it is even only necessary that $\|\varepsilon_0\|_{L^{2}}=o(\delta x^{1/2})$.
\begin{flushright}
$\Box$
\end{flushright}

\begin{theorem}
\label{thstabBurgers}
If we solve the Burgers equation (\ref{burgers_eqn}) with the numerical scheme (\ref{RKtype}), then for a sufficiently regular solution $u$,
the stability condition is provided by:
\begin{equation}
\label{CFLthBur}
\delta t\leq \left( \frac{2(C-\|\partial_x u\|_{L^\infty})}{S_r}\right)^{\frac{1}{2r-1}} \left( \frac{\delta x}{\|u\|_{L^{\infty}}} \right)^{\frac{2r}{2r-1}}.
\end{equation}
where $\delta t$ is the time step, $\delta x$ the space step, $r$ an integer and $S_r$ a quantity both defined by
Eq. (\ref{Aoxi}), (\ref{sumS_q}) and (\ref{calcSell}) and
$C$ the constant in the exponential growth of the error: $\varepsilon(t)\sim\varepsilon_0e^{Ct}$.
\end{theorem}
\emph{proof:} Thanks to proposition \ref{thRn}, we are able to write:
\begin{equation}
\|\varepsilon_{n+1}\|_{L^2}
\leq\|\varepsilon_{n}+\sum_{i=1}^s \beta_{i} \delta t^i u_n^i\partial_x^i\varepsilon_{n}\|_{L^2}
+\delta t \|\varepsilon_n\partial_x u_n\|_{L^2}+o(\delta t \|\varepsilon_n\|_{L^2})
\end{equation}
On the other hand we have $\|\varepsilon_n\partial_x u_n\|_{L^2}\leq \|\partial_x u\|_{L^\infty}\|\varepsilon_n\|_{L^2}$ and
since for $i,j\leq s$,
\begin{equation}
\label{BurFF}
\langle\delta t^i u_n^i\partial_x^i\varepsilon_{n},\delta t^j u_n^j\partial_x^j\varepsilon_{n}\rangle_{L^2}=\left\{
\begin{array}{ll} o(\delta t \|\varepsilon_n\|_{L^2}^2) &~~~~{\rm if}~~~~i+j=2\ell+1 \\ & \\
(-1)^{\ell-i}\delta t^{2\ell}\|u_n^\ell\partial_x^\ell\varepsilon_{n}\|_{L^2}^2+o(\delta t \|\varepsilon_n\|_{L^2}^2) &~~~~{\rm if}~~~~i+j=2\ell
\end{array}\right.
\end{equation}
we derive
\begin{equation}
\|\varepsilon_{n}+\sum_{i=1}^s \beta_{i} \delta t^i u_n^i\partial_x^i\varepsilon_{n}\|_{L^2}^2
=\sum_{\ell=0}^{2s} S_\ell \delta t^\ell \|u_n^\ell\partial_x^\ell\varepsilon_{n}\|_{L^2}^2+o(\delta t \|\varepsilon_n\|_{L^2}^2)
\end{equation}
with $S_\ell$ given by (\ref{calcSell}).\\
Then, as $S_0=1$, and
$\|u_n^\ell\partial_x^\ell\varepsilon_{n}\|_{L^2}\leq\|u_n\|_{L^{\infty}}^\ell\frac{\|\varepsilon_{n}\|_{L^2}}{\delta x^\ell}$,
\begin{equation}
\|\varepsilon_{n}+\sum_{i=1}^s \beta_{i} \delta t^i u_n^i\partial_x^i\varepsilon_{n}\|_{L^2}^2
\leq\sum_{\ell=0}^{s} S_\ell \left(\frac{\delta t}{\delta x}\right)^{2\ell}
\|u\|_{L^{\infty}}^{2\ell}\|\varepsilon_{n}\|_{L^2}^2+o(\delta t \|\varepsilon_n\|_{L^2}^2)\\
\end{equation}
so
\begin{equation}
\|\varepsilon_{n}+\sum_{i=1}^s \beta_{i} \delta t^i u_n^i\partial_x^i\varepsilon_{n}\|_{L^2}
\leq\left(1+\frac{1}{2}\sum_{\ell=1}^{s} S_\ell \left(\frac{\delta t}{\delta x}\right)^{2\ell}
\|u\|_{L^{\infty}}^{2\ell}+o(\delta t)\right)\|\varepsilon_{n}\|_{L^2}
\end{equation}
and finally
\begin{equation}
\|\varepsilon_{n+1}\|_{L^2}
\leq\left(1+\frac{1}{2}\sum_{\ell=1}^{s} S_\ell \left(\frac{\delta t}{\delta x}\right)^{2\ell}
\|u\|_{L^{\infty}}^{2\ell}+\delta t \|\partial_x u\|_{L^\infty}+o(\delta t)\right)\|\varepsilon_{n}\|_{L^2}.
\end{equation}
Let $r$ be the first power in the sum where $S_r\neq 0$, and let us assume that $S_r>0$.
Then, the stability condition $\|\varepsilon_{n+1}\|_{L^2}\leq\left(1+C\delta t\right)\|\varepsilon_{n}\|_{L^2}$
is reduced to
\begin{equation}
\frac{1}{2} S_r \frac{\delta t^{2r-1}}{\delta x^{2r}}\|u\|_{L^{\infty}}^{2r}\leq (C-\|\partial_x u\|_{L^\infty})
\end{equation}
i.e. the condition (\ref{CFLthBur}).
We recognize the same power law as the one obtained in the linear case (\ref{CFLdtdx}).
The term $-\|\partial_x u\|_{L^\infty}$ should usually be discarded since its contribution is external
to the instability phenomenon and random.

\subsection{Scalar conservation laws}
\label{cons-law}
Scalar conservation laws group equations of the type
\begin{eqnarray}
\label{scal_laws}
\partial_t u + \sum_{i=1}^d \partial_{x_i} f_i(u) = 0\quad {\rm for}~~(\x,t)\in\R^d\times[0,T]\\
u(0,\x)=u_0(\x)\quad {\rm for}~~\x\in\R^d
\end{eqnarray}
with $f_i:\R\to\R$ differentiable functions and $u:\R^d\to\R$ the scalar unknown function.\\
A stability analysis of the solution of these equations in the frame of Discontinuous Galerkin Runge-Kutta
formulation was presented in \cite{ZS04} for space accuracy of order two and three, with the $\delta t\leq C\delta x^{4/3}$
CFL-like condition, but as the byproduct of a long and rigorous computational process.
This work was the continuation of \cite{CS01} where the authors observed that first and second order Runge-Kutta
methods are unstable under any linear CFL conditions when the space discretization is sufficiently accurate
and so does not dissipate too much. In this section, we link their results to our analysis and refine
the stability criteria. Actually we have the following result:

\begin{theorem}
\label{thscl}
Let us apply the numerical scheme (\ref{RKtype}) to solve the equation (\ref{scal_laws}).
If $f\in C^{p+1}$ and $u\in C^{p}$ i.e. $f^{(p+1)},u^{(p)}\in L^\infty$ and if the $(S_\ell)$ defined by
(\ref{Sell}) satisfy $S_1=\dots=S_{r-1}=0$ and $S_r>0$,
then, given a constant $C$ limiting the exponential growth of the stability error: $\varepsilon_T\leq e^{CT}\varepsilon_0$,
the numerical scheme is conditionally stable under the CFL-like condition:
\begin{equation}
\label{eqthscl}
\delta t\leq \left( \frac{2C}{S_r}\right)^{1/(2r-1)} \left( \frac{\delta x}{\sum_{i=1}^d\|f_i'(u)\|_{L^{\infty}}} \right)^{\frac{2r}{2r-1}}.
\end{equation}
\end{theorem}
{\sl proof:}
The proof is more or less the same as for the Burgers case {\it cf} part \ref{burger}, using the following facts:
\begin{itemize}
\item $f_i(u_n+\varepsilon_n)=f_i(u_n)+f_i'(u_n)\varepsilon_n+o(\varepsilon_n)$,
\item $u_{(\ell)}-u_n=o(1)$,
\item $\partial_{x_i}\left( f_i'(u_n)\varepsilon \right)\sim f_i'(u_n)\partial_{x_i}\varepsilon$ for stability analysis, and
\item for all functions $\eta$ and $\varepsilon$,
\begin{equation}
\left\langle \eta,\sum_{i=1}^d\sum_{j=1}^d \partial_{x_i}\left( f_i'(u_n)f_j'(u_n)\partial_{x_j}\varepsilon \right) \right\rangle
=-\left\langle \sum_{i=1}^d f_i'(u_n)\partial_{x_i}\eta,\sum_{i=1}^d f_i'(u_n)\partial_{x_i} \varepsilon \right\rangle
\end{equation}
allowing equalities of the type (\ref{FF}).
\end{itemize}
Finally, we obtain:
\begin{equation}
\|\varepsilon_{n+1}\|_{L^2}^2=(1+2C_1\delta t+o(\delta t))\|\varepsilon_n\|_{L^2}^2
+{S_r}\,{\delta t^{2r}}\left\|\sum_{{\bf i}\in[1,d]^r}\left( \prod_{s=1}^rf_{i_s}'(u_n)\right)
\left( \prod_{s=1}^r\partial_{x_{i_s}}\right) \varepsilon_n\right\|_{L^2}^2
\end{equation}
Then, knowing that for $\varepsilon_n\in V_{\delta x}$,
\begin{eqnarray}
\nonumber\left\|\sum_{{\bf i}\in[1,d]^r}\left( \prod_{s=1}^rf_{i_s}'(u_n)\right)
\left( \prod_{s=1}^r\partial_{x_{i_s}}\right) \varepsilon_n\right\|_{L^2}^2
&\leq &\left( \sum_{{\bf i}\in[1,d]^r} \left( \prod_{s=1}^r\|f_{i_s}'(u_n)\|_{L^\infty}\right)
\frac{\|\varepsilon_n\|_{L^2}}{\delta x^r} \right)^2\\
&\leq &\left( \left( \sum_{i\in[1,d]}\|f_{i}'(u_n)\|_{L^\infty} \right)^r \frac{\|\varepsilon_n\|_{L^2}}{\delta x^r} \right)^2,
\end{eqnarray}
and neglecting the constant $C_1$, the von Neumann stability criteria
\begin{equation}
\|\varepsilon_{n+1}\|_{L^2}^2\leq(1+2C\delta t+o(\delta t))\|\varepsilon_n\|_{L^2}^2
\end{equation}
is satisfied if
\begin{equation}
\left( \sum_{i\in[1,d]}\|f_{i}'(u)\|_{L^\infty} \right)^{2r} \frac{S_r\,\delta t^{2r}}{\delta x^{2r}}
\leq 2C\delta t,
\end{equation}
i.e. condition (\ref{eqthscl}).

\subsection{Incompressible Euler equation}
\label{euleri}
The Euler equations model incompressible fluid flows with no viscous term:
\begin{eqnarray}
\label{Eulervp}
\frac{\partial \mathbf{u}}{\partial t}+({\bf u}\cdot\nabla) {\bf u}-\nabla p=0,\quad 
{\rm div}\,\u=0,\quad {\rm for}~(t,\x)\in\R_+\times \Omega.
\end{eqnarray}
The use of the Leray projector $\mathbb{P}$ which is the $L^2$-orthogonal projector on the divergence-free
space, allows us to remove the pressure term:
\begin{eqnarray}
\label{EulervPdiv}
\frac{\partial \mathbf{u}}{\partial t}+\mathbb{P}\left[({\bf u}\cdot\nabla) {\bf u}\right]=0.
\end{eqnarray}
The stability analysis of this case proceeds somehow as a synthesis of the previous two sections Sec. \ref{transpU} and Sec. \ref{burger}.
An important property is then the skewness property of the transport term
(see \cite{GR86}, chapter IV, Lemma 2.1 or \cite{DP07} for the proof, also used for the stability of
the incompressible Navier-Stokes equations in \cite{MT98}):
\begin{lemma}
\label{lemme}
Let $\u,{\bf v},\w\in H^1(\Omega)^d$, $H^1(\Omega)$ denoting the Sobolev space on the open set $\Omega\subset \R^d$,
be such that $(\u\cdot\nabla){\bf v},(\u\cdot \nabla)\w \in L^2$.
If $\u\in {\textbf{H}}_{{\rm div},0}(\Omega)=\{\f\in (L^2(\Omega))^d,~{\rm div}~\f=0~{\rm on}~\Omega,~\f\cdot \n=0~{\rm on}~\partial \Omega\}$, then
\begin{equation}
\label{skewsymrel}
\langle {\bf v},(\u\cdot \nabla)\w\rangle_{L^2(\Omega)}=-\langle(\u\cdot \nabla){\bf v},\w\rangle_{L^2(\Omega)}.
\end{equation}
\end{lemma}

\begin{corollary}
\label{corollary}
With the same assumptions as in lemma \ref{lemme},
\begin{equation}
\langle{\bf v},(\u\cdot \nabla){\bf v}\rangle_{L^2(\Omega)}=\int_{\x\in\Omega}{\bf v}\cdot(\u\cdot \nabla){\bf v}\, d\x=0.
\end{equation}
\end{corollary}

\ 

Considering the scheme (\ref{RKtype}), we introduce a stability error $\varepsilon_{(\ell)}$ at level $\ell$.
Then, under the condition $\delta t=o(\delta x)$ and for $\varepsilon_{n}$ small enough,
most of the terms appearing in the expression of $\varepsilon_{(\ell)}$ are
negligible with respect to:
\begin{itemize}
\item the terms $\delta t^i\,F_n^i(\varepsilon_{n})$ where
$F_n(\varepsilon_n)=\P[(\u_n\cdot\nabla)\varepsilon_n]$ and
$F_n^i=\underbrace{F_n\circ F_n\circ \dots \circ F_n}_{i~{\rm times}}$,
\item the term $\delta t\, \P[(\varepsilon_{n}\cdot\nabla)\u_n]$.
\end{itemize}
Then most of the arguments used in Sec. \ref{burger} apply with even more accuracy since
we have the orthogonality relation:
\begin{equation}
\label{FF}
\langle F_n^i(\varepsilon_n),F_n^j(\varepsilon_n)\rangle_{L^2(\Omega)}=\left\{
\begin{array}{ll} 0 &~~~~{\rm if}~~~~i+j=2\ell+1~~~~{\rm for}~~~~\ell\in\N \\
(-1)^{\ell-i}\|F_n^\ell(\varepsilon_n)\|_{L^2(\Omega)}^2 &~~~~{\rm if}~~~~i+j=2\ell~~~~{\rm for}~~~~\ell\in\N
\end{array}\right.
\end{equation}
instead of (\ref{BurFF}).
This leads to the following result:
\begin{proposition}
\label{stabEul}
Assume that the incompressible Euler equations (\ref{Eulervp}) have a $s$-times space-differentiable
solution $\u$ such that $\|\nabla^s \u\|_{L^{\infty}([0,T]\times\Omega)}<+\infty$,
that the discretization conserves the skew-symmetry relation (\ref{skewsymrel})
and that $\forall\varepsilon\in V_{\delta x}$,
$\|\P_{\delta x}\,F_n^k(\varepsilon)\|_{L^2}\leq C\|\u\|_{L^{\infty}}^k \frac{\|\varepsilon\|_{L^2}}{\delta x^k}$.
Then a stability error $\varepsilon$ small enough at the initial time: $\|\varepsilon_0\|_{L^{2}}=o(\delta x^{d/2})$
remains bounded for $t\in[0,T]$ under the condition:
\begin{equation}
\label{stabEuleq}
\delta t\leq \left( \frac{2C}{S_r}\right)^{1/(2r-1)} \left( \frac{\delta x}{\|\u\|_{L^{\infty}}} \right)^{\frac{2r}{2r-1}}
\end{equation}
with $\delta t$ the time step, $\delta x$ the space step, and $r$ and $S_r$ obtained as in (\ref{Sell}).
\end{proposition}

\ 

This proposition extends to Navier-Stokes equations for high Reynolds number.
The incompressible Navier-Stokes equations are written:
\begin{equation}
\label{NS}
\left\{
\begin{array}{l}
\partial_t {\bf u} + {\bf u}\cdot \nabla {\bf u}-\nu \Delta {\bf u}+\nabla p=0,\\
{\rm div} {\bf u}=0,\\
{\bf u}(0,x)={\bf u}_0(x)
\end{array}
\right.
~x\in\R^d,~t\in[0,T]
\end{equation}
Using the Leray projector $\mathbb{P}$ --the orthogonal projector on divergence-free vector fields--
we reduce the equation to:
\begin{equation}
\label{NSwP}
\partial_t\mathbf{u}+\mathbb{P}\left[\mathbf{u}\cdot
\nabla\mathbf{u}\right]-\nu\Delta\mathbf{u}=0
\end{equation}
Two second order schemes are widely in use for the solution of this equation: the order two Runge-Kutta scheme \cite{KT97,DP07}
as well as the second order Adams-Bashforth scheme \cite{Pey00,Sch05}.

When the Reynolds number $Re=\frac{\|u\|_{L^\infty}L}{\nu}$ is sufficiently large, the contribution of the heat kernel to
the stability vanishes \cite{DP07}, and the same instability effects as for the incompressible Euler
equation appear as it was observed in 2D experiments \cite{DP07}.
New tests with boundaries comply the stability condition $\delta t\leq\delta t_{\max}= C \delta x^{4/3}$ for dipole/wall numerical
experiments.
These results will be presented in a forthcoming paper.

\section{Multi-component transport}
\label{crossing}
We extend the scope of application of the stability conditions (\ref{CFLdtdx}) to other cases with multiple
derivatives in time, like wave equations, or multiple components, like in some MHD models \cite{CH08}.
Let us consider the one dimensional equation:
\begin{equation}
 \partial_t X=M \partial_x X,\quad {\rm with}~~X=\left( \begin{array}{c} u_1 \\ u_{2} \\ \vdots \\ u_{n} \end{array} \right),
\quad u_\ell:\R\to \R,\quad {\rm and}~~ M\in \mathcal{M}_n(\R).
\end{equation}
For example, for $X=(u,v)^t$, and $M=\left[ \begin{array}{cc} 0 & 1 \\ 1 & 0 \end{array} \right]$, we obtain the
wave equation $\partial_t^2 u=\partial_x^2 u$.\\
Regarding the general case, we diagonalize the matrix $M$ in $\C$:
\begin{equation}
M=P^{-1}DP,\quad {\rm with}~~D=\left[ \begin{array}{ccc} \lambda_1 & & 0\\  & \ddots &  \\ 0 & & \lambda_n \end{array} \right].
\end{equation}
Considering $Y=PX$, the equation $\partial_t Y=D \partial_x Y$ has physical meaning only if $\lambda_\ell\in\R$ for all $\ell$.
Under this form all the components are independent.
Therefore all our results on the transport equation apply to this case taking $a=\max_\ell|\lambda_\ell|$.

\ 

When the matrix $M$ cannot be diagonalized, like in the case
$M=\left[ \begin{array}{cc} 1 & 1 \\ 0 & 1 \end{array} \right]$,
we remark that the second component is independent from the first component:
\begin{equation}
\left\{ \begin{array}{l} \partial_t u_1= \partial_x u_1 + \partial_x u_2 \\ \partial_t u_2= \partial_x u_2 \end{array} \right.
\end{equation}
then the term $\partial_x u_2$ in the first equation plays the role of a source term.

\ 

In the case when there are several space variables:
\begin{equation}
 \partial_t X=M_1 \partial_{x_1} X+M_2 \partial_{x_2} X+\dots+M_n \partial_{x_n} X
\end{equation}
applying a von Neumann stability analysis, we obtain:
\begin{equation}
 \partial_t \hat X=\left( M_1\,i\xi_1+M_2\,i\xi_2 +\dots+M_n\,i\xi_n \right) \hat X.
\end{equation}
We consider $M(\xi)=M_1\,\xi_1+M_2\,\xi_2 +\dots+M_n\,\xi_n$, and diagonalize $M(\xi)=P(\xi)^{-1}D(\xi)P(\xi)$.
As previously, taking $\hat Y_\xi=P(\xi)\hat X$, we obtain the stability constraint (\ref{eqFLST}) as in the scalar case.

\section{Conclusion}
\label{concl}
The stability CFL-like conditions presented in this paper may be encountered in many simulations
of convection-dominated problems using explicit numerical schemes.
Although based on a classical von Neumann stability analysis,
this kind of stability analysis is not performed elsewhere.

Two arguments support our approach.
First we explain some ``CFL shrinking'' effects for second order schemes
already in use: people remarked that they had to take $C\to 0$ in the usual linear CFL
condition $\delta t\leq C\delta x$.
Secondly, we predict some exotic CFL conditions $\delta t\leq C\delta x^{\frac{2r}{2r-1}}$ for
certain Runge-Kutta schemes and Adams-Bashforth schemes which optimize the energy conservation. 
Numerical tests validate these predictions.

We showed why increasing the temporal order of a scheme increases the stability.
We even linked the order of a scheme, its stability and the tangency of its stability
domain to the $(Oy)$ axis in the von Neumann stability analysis.
Nevertheless, the numerical viscosity may erase these instability effects especially when
using an upwind scheme \cite{CS01}.

We extended the domain of application of these results to different equations, including equations on bounded domains,
non linear equations, and equations with multiple derivatives in time.
These extensions assume smoothness properties for the solution.
This smoothness assumption restrains the frame of application to a rather limited area.
Nevertheless, this clear exposing of actual numerical artifacts plus the correlate accurate stability conditions should
be useful to a wide community, in particular
to those who perform numerical simulations of turbulent flows with spectral codes.

\section*{Acknowledgements}

The author gratefully acknowledges the CEMRACS 2007 organizers for his stay in the CIRM in Marseilles
and for his access to its rich bibliographical resources, as well as Institute of Fundamental Technological
Research Polish Academy of Sciences (IPPT PAN) and Commissariat \`a l'\'Energie Atomique (CEA) for his stays there
in 2007/2008 and 2009 respectively.
He wishes to express his gratitude to Yvon Maday and Fr\'ed\'eric Coquel
for fruitful discussions, as well as to Dmitry Kolomenskiy for his help in the redaction of this paper and
in the realization of the numerical experiments.
He also acknowledges the anonymous referees, whose comments substantially improved the quality of the paper.


\bibliographystyle{plain}

 

 
\end{document}